\documentclass{amsart}
\usepackage{xypic,amsmath,amssymb}
\def\demo{\par\noindent{\bf Proof.\ }}
\def\enddemo{\ $\Box$\par\vskip.6truecm}
\newtheorem{definition}{Definition}[section]
\newtheorem{theorem}{Theorem}[section]
\newtheorem{proposition}{Proposition}[section]
\newtheorem{oss}{Remark}[section]
\newtheorem{corollary}{Corollary}[section]
\newtheorem{lemma}{Lemma}[section]

\newtheorem{esempi}{Examples}[section]

\newcommand{\homs}{{\sf Hom}}

\newcommand{\hsolo}{\mathcal{H}_s}
\newcommand{\holo}{\mathcal{H}}
\newcommand{\hsolop}{{\mathcal{H}_s}_\bullet}
\newcommand{\holop}{{\mathcal{H}}_\bullet}
\newcommand{\fasc}{\boldsymbol{\mathcal{F}}_T(\mathcal S)}
\newcommand{\fasci}{\boldsymbol{\mathcal{F}}(\co)}

\newcommand{\ip}{\mathfrak{Ip}}

\newcommand{\co}{\mathfrak{Compl}}

\newcommand{\tlc}{Tlc_{open}}

\newcommand{\pro}{\mathbb{P}^1}
\newcommand{\ycal}{\mathcal{Y}}
\newcommand{\xcal}{\mathcal{X}}
\newcommand{\zcal}{\mathcal{Z}}
\newcommand{\fcal}{\mathcal{F}}

\newcommand{\ccal}{\mathcal{C}}
\newcommand{\scal}{\mathcal{S}}
\newcommand{\aff}{\mathbb{A}^1}
\newcommand{\yon}{\mathbf{Y}}

\newcommand{\compl}{\mathbb{C}}

\newcommand{\proj}{\mathbb{P}}

\newcommand{\ins}{{\sf Ins}}
\newcommand{\topo}{{\sf Top}}

\newcommand{\fs}{\mathfrak{S}}

\newcommand{\pts}{\mathsf{pt}}
\begin{document}
\baselineskip=.7cm
\title[Extended hyperbolicity]{Extended hyperbolicity}
%\author[Borghesi and Tomassini]{Simone Borghesi$^{*}$ \and
%Giuseppe Tomassini$^{*}$}\thanks{$^{*}$Supported by the MURST
%project \lq\lq Geometric Properties of Real and Complex Manifolds".}
\author[Borghesi and Tomassini]{Simone Borghesi\and Giuseppe
  Tomassini}
\thanks{Supported by the MURST project \lq\lq Geometric Properties of
Real and Complex Manifolds".}


%\sc Scuola Normale Superiore di Pisa, Piazza dei Cavalieri

%7, 56126 Pisa, Italy}
\maketitle
%\centerline{\em 1st of June, 2007}
\begin{abstract}
Given a complex space $X$, we cosidered the problem of finding a 
{\it hyperbolic model} of $X$. This is an
object $\ip(X)$ with a morphism $i:X\to \ip(X)$ in such a way
that $\ip(X)$ is ``hyperbolic'' in a suitable sense and $i$ is as
close as possible to be an isomorphism. Using the theory of model
categories, we found a definition of hyperbolic simplicial sheaf (for
the strong topology) that extends the classical one of Brody for complex spaces. We prove the existence of hyperbolic models for any simplicial
sheaf. Furthermore, the morphism $i$ can be taken to be a cofibration
and an affine weak equivalence (in an algebraic setting, Morel and
Voevodsky called it an $\aff$ weak equivalence). Imitating one
possible definition of homotopy groups for a topological space, we
defined the {\it holotopy} groups for a simplicial sheaf and showed
that their vanishing in ``positive'' degrees is a necessary condition for
a sheaf to be hyperbolic.  We deduce that if $X$ is a complex space
with a non zero holotopy group in positive degree, then its hyperbolic
model (that in general will only be a simplicial sheaf) cannot be
weakly equivalent to a hyperbolic complex space (in particular is not
itself hyperbolic). We finish the manuscript by applying these results
and a {\it topological realization functor}, constructed in the
previous section, to prove that the
hyperbolic models of the complex projective spaces cannot be weakly
equivalent to  hyperbolic complex spaces.
\end{abstract}

%%%%%%%%%%%%%%%%%%%%%%%%%%%%%%%%%%%%%%%%%%%%%%%%%%
%%%%%%%%%%%%%%%%%%%%%%%%%%%%%%%%%%%%%%%%%%%%%%%%%%%

\tableofcontents

\section{Introduction}
The notion of hyperbolic space was given by Kobayashi in \cite{koba1}. It is 
based on the existence of certain intrinsic distances, originally introduced 
to generalize  Schwarz Lemma to higher dimensional complex spaces.
Let $D\subset \mathbb C$ be the unit disc endowed with the Poincar\'e distance 
$\rho$. In view of the Schwarz-Pick Lemma every holomorphic map $f:D\to D$ is a 
contraction for $\rho$. Let $X$ be a complex space.  A {\it chain of holomorphic 
discs} between two points $p$,$q$ in $X$ is a set $\alpha=\{f_1,\cdots, f_k\}$ 
of holomorphic maps $D \to X$  such
 that there are points $p=p_0, p_1,\cdots, p_k=q$ in $X$ and $a_1,b_1, \cdots,
  a_k, b_k$ in $D$ with the property $f_i(a_i)=p_{i-1}$ and $f_i(b_i)=p_i$,
  $i=1,\cdots, k$. The {\it length} of $\alpha$ is defined as
\begin{equation}
l(\alpha)=\sum_{i=1}^k\rho(a_i,b_i)
\end{equation}
The {\it Kobayashi pseudo distance} $d_{Kob}$ on $X$ is defined as
\begin{equation}
d^X_{Kob}(p,q)=\inf_{\alpha}~ l(\alpha)
\end{equation}
where $\alpha$ varies through the family of all chains of holomorphic discs
joining $p$ and $q$. For quasi-projective varieties the pseudo distance of Kobayashi can be defined 
by means of chains of algebraic curves (see \cite{dls}).

The contraction property holds with respect to the Kobayashi pseudodistances: if  $f: X\to Y$ is a holomorphic map
between complex spaces, we have 
\begin{equation}\label{contr}
d_{Kob}^Y(f(p),f(q))\leq d_{Kob}^X(p,q),
\end{equation}
for every $p,q\in X.$ In particular, $d_{Kob}^Y$ is invariant by 
biholomorphisms. It follows that $d^D_{Kob}=\rho$.

We have $d_{Kob}^{\compl}\equiv 0$, $d_{Kob}^{\compl^\ast}\equiv 0$. More generally, 
one has $d_{Kob}^G\equiv 0$ for every connected
complex Lie group $G$ (see \cite {koba}). 

A complex space $X$ is said to be {\em hyperbolic} (in the sense of
Kobayashi) if $d^X_{Kob}$ is a distance. The unit disc $D$ is
hyperbolic, whereas $\compl$ is not. $\compl\setminus A$ with
${\sf card}\,A\ge 2$ is hyperbolic. A compact complex curve of genus 
$g\ge 2$ is a hyperbolic space \cite{koba}.
 $X$ is said to be {\em hyperbolic modulo C}, where $C$ is a closed 
 subset (usually a closed complex 
subspace), if for every pair of distinct points $x, y\in X\setminus  C$ we have $d^X_{Kob}(x,y)>0$. 

If a complex space $Y$ is $\compl$-connected (i.e. for any
$p\neq q$ points in $Y$ there exists a holomorphic function
$f:\compl\to Y$ such that $p,\,q,\in f(\compl)$) then, by virtue
of the contraction property the only holomorphic maps with values in a hyperbolic 
space $X$ are the constant ones. In particular, every holomorphic map $\mathbb{C}\to X$ 
is constant. The crucial fact is that for compact complex spaces the converse is also 
true. This is the content of the fundamental theorem of Brody (cfr.   \cite{koba}, \cite{lang}).

This result motivates the following definition: a complex space $X$ is said to be 
{\em Brody hyperbolic} if every holomorphic map $f:\compl\to X$ is constant.

As well as for hyperbolicity we have the notion of {\it Brody 
hyperbolicity modulo a closed subset $C$}: $X$ is said to be
{\em Brody hyperbolic modulo C} if every non constant holomorphic map
$f:\compl\to X$ satisfies $f(\compl)\subset C$.
A Kobayashi hyperbolic space is Brody hyperbolic but the converse is in general 
not true. Indeed Mark Green constructed a Zariski open set $W$ in $\mathbb{P}^2$ 
(the two dimensional complex projective space) , 
deleting four lines in general position and three points outside the four lines, 
which is Brody hyperbolic but not Kobayashi hyperbolic \cite{lang}.

Related to hyperbolicity are some basic conjectures which motivated several important 
papers in Algebraic and Analytic Geometry. \begin{enumerate}
\item[(1)] A generic hypersurface of degree $\ge 2n+1$ in $\proj^n$ is 
hyperbolic;
\item[(2)] The complement of a hypersurface of degree $\le 2n$ in 
$\proj^n$ is not hyperbolic;
\item[(3)] A generic hypersurface of degree $\ge n+2$ in $\proj^n$ 
is hyperbolic modulo a proper closed subvariety;
\item[(4)] a smooth projective hyperbolic variety has an ample canonical bundle 
(Kobayashi's conjecture);
\item[(5)] a smooth algebraic variety is of general type if and only if it is 
hyperbolic modulo a proper algebraic subset (Lang's conjecture \cite{lang1}, \cite{lang2}). 
\end{enumerate}
For the basic material as well as a discussion of 
the geometric meaning of these conjectures we refer to \cite{koba}, \cite{lang},  
\cite{lang1}, \cite{lang2} and the bibliography there.

In this paper we will consider the following problem: given a complex
space $X$, construct a ''hyperbolic model" of $X$ i.e. a ``hyperbolic'' object 
$\mathfrak{Ip}(X)$, in a sense the ``closest'' hyperbolic object to $X$ endowed 
with a canonical natural map $c_X:X\to \mathfrak{Ip}(X)$ having the following 
universal property: if $Y$ is hyperbolic a holomorphic map $f:X\to Y$ factorizes 
through $\mathfrak{Ip}(X)$ i.e. we have a commutativity diagram 
$$
\xymatrix{X\ar[r]^{c_X}\ar[d]^f &
  \mathfrak{Ip}(X)\ar[dl]^{\tilde{f}} \\
Y & }.
$$
One possible way to do this would be to consider the quotient topological space
$X/\mathcal R$ where $\mathcal R$ is the equivalence relation: $x\sim y$ iff $
d_{Kob}^X(x,y)=0$ or, bearing in mind Brody's Theorem, if and only if they belong 
to the image of a holomorphic map $\compl\to X$. 
This approach has two oddnesses. One is that 
$X/\mathcal R$ is in general very different from $X$, indeed $X/\mathcal R$
is just a point for $\compl$-connected spaces $X$. On the other hand, $X/\mathcal R$ 
will have in general no complex structure (even in a weak sense), thus it will 
be impossible to define a Kobayashi psedodistance on this quotient in order to have an 
useful concept of hyperbolicity on it. 

Regarding this, it is worth mentioning the nice paper of Campana \cite{camp} where a concept of Kobayashi pseudodistance is defined for orbifolds.
Then, for any variety which is smooth  
and bimeromorphically equivalent to a  K\"ahler manifold he constructed an orbifold $C(X)$ called the {\it core} of $X$ and a meromorphic function $c_X: X\to C(X)$. 
Furthermore, he conjectured that the 
generic fiber of $c_X$ has a vanishing Kobayashi metric and $C(X)$ is Brody 
hyperbolic modulo a proper subvariety.

In this paper we developed a different approach. We used techniques pioneered by Quillen in
\cite{quillen} and largely employed in \cite{morvoe}, which we drew
inspiration from in writing the technical sections of this paper. 
We construct an (unstable) homotopy category of complex spaces $\holo$, whose objects
include (homotopy) classes of complex spaces. Unlike the classical homotopy category
of topological spaces, the category
$\holo$ reflects the complex structure of the objects. The procedure involves an
enlargement of the category of complex spaces to a new category containing as
full subcategory the one of complex spaces with holomorphic functions. In this 
bigger category we define a notion of hyperbolicity which we prove that
it restricts to the Brody hyperbolicity for complex spaces. 
Using this notion, we show that in each class of complex spaces lies a hyperbolic
representative $\mathfrak{Ip}(X)$, which in general will not be a complex space. 
It follows that $\mathfrak{Ip}$ will be (weakly equivalent to) a point
if and only if $X$ is. 
Such correspondence is functorial and there exists a canonical morphism
$c_X:X\to \mathfrak{Ip}(X)$ satisfying the universality property described
above. $c_X$ and $\widetilde{f}$ will be morphisms of the homotopy category
in general, but the composition $\widetilde{f}\circ c_X$ is a class represented
by a holomorphic function and the commutativity is as holomorphic functions
as opposed to "`up to homotopy"'.
Concerning the object $\mathfrak{Ip}(X)$, we will prove that the class of $\mathbb{P}^n$ 
cannot have a hyperbolic complex space as representative, whereas in the class of
$\compl$, the point can be taken as hyperbolic complex space
representative. $\mathfrak{Ip}(X)$ is given by a complicated
construction even if $X$ is a complex space, although its {\it topological realization} 
(see Section 6) is a topological space homotopic equivalent to the topological
space underlying $X$. 

The procedure to construct the category $\holo$ follows closely the
one described in \cite{morvoe} which works in a quite general
context. It follows that almost all the results proved here are valid for algebraic schemes of 
finite type over a noetherian base, as well.

The main idea is to construct a category
obtained from another by ``adding'' the inverses of certain morphisms.
In the case of the category $\co$ of complex spaces with the strongly topology and holomorphic
maps, we wish to add the inverse to the canonical map
$p:\compl\to {\sf pt}$ (the canonical projection $\aff_B\to B$ in the
algebraic case) along with all its base changed maps. Such a category, 
which we denote as
$p^{-1}\co$, exists, however, to make it usable, it should be
obtained as the homotopy category associated to a {\em model
  structure} (see \cite{quillen}) on $\co$. In general, deciding whether
a localized category $S^{-1}\mathcal{C}$ is equivalent to the homotopy
category associated to a model structure on $\mathcal{C}$ is a very
complicated task. This has been proved in the case of derived
categories and the homotopy category of topological spaces. There are
only partial results on this issue, if we assume that the category
$\mathcal{C}$ is a homotopy category itself and posseses a ``simplicial
structure''. The easiest way to replace a category
$\mathcal{C}$ with one endowed of such simplicial structure is to
consider $\Delta^{op}\mathcal{C}$, the category of simplicial objects
in $\mathcal{C}$. Then, we may try to give to $\Delta^{op}\mathcal{C}$
a {\em simplicial} model structure. If all this is successful, the
homotopy category associated to such simplicial model structure is a
good candidate to start with for establishing whether we can localize
with respect of some morphism by using an appropriate model structure.
In our situation, the category $\co$ is replaced with $\fasc$, the
category of sheaves over the site of complex spaces endowed with the
strong topology $T$ (in the algebraic case, this will denote the category
of sheaves over the site of smooth schemes of finite type over a
noetherian base endowed with a topology not finer than quasi compact
flat topology). The reason for doing this lies mainly in the fact that
not all diagrams admit colimits and
the existing ones in $\co$ often are unsuitable to do homotopy
theory with (see Section \ref{colimiti} for more details on this). 
On the other hand, $\fasc$ is complete and cocomplete
and the colimits have a ``suitable'' shape. We than proceed with the
program described above in order to invert $p:\compl\to \sf{pt}$.
We end up with the category $\hsolo$ which is defined as the homotopy
category associated to the simplicial model structure on
$\Delta^{op}\fasc$. The morphism $p$ in the category $\hsolo$ fits in
the Bousfield framework \cite{bousfield}, and lies inside the class of weak
equivalences in an appropriate model structure on
$\Delta^{op}\fasc$. The associated homotopy category will be denoted by
$\holo$ and sometimes by $\holo^{olo}$ when we wish to stress that
we are in the holomorphic setting. Any object of the site represents a
class in $\holo$ and in the case it is a complex space, its hyperbolic
model $\widetilde{X}$ will be only a simplicial sheaf on
$\co$. The notion of hyperbolicity for a simplicial sheaf $\xcal$ is given in
the Definition \ref{spazio-iper}. In the particular case $\xcal=X$ is
a compact complex space, in view of Theorem \ref{fibering} and 
Brody's Theorem, we conclude that $X$ is hyperbolic according to our
definition if and only if it is Kobayashi hyperbolic (see Corollary
\ref{equivalenti}). Thus, the Definition \ref{spazio-iper} is a
generalization of the classical concept of hyperbolicity for complex
spaces.

In the  section \ref{Holotopy/G} we associate certain sets to each
object of $\holo^{olo}$
which have a natural group structure in positive simplicial
degree. They are  called {\em holotopy sets} or {\em groups} when
applicable (see Definition \ref{gruppi}). We prove that the vanishing
of some of the holotopy groups of a complex space $X$ is a necessary
condition for the hyperbolic model $\ip(X)$ to be isomorphic in
$\holo^{olo}$ to some hyperbolic complex space.

In the following section we construct an useful functor for explicit 
computations: the {\it topological realization functor}. To a
simplicial sheaf it associates a  topological space in such a way 
few reasonable properties are satisfied (cfr. Definition 
\ref{real-topo}). In the last section, as an application of some of
our results, we show that 
$\ip({\mathbb{P}^n})$ is not weakly equivalent to a Brody hyperbolic 
space for any $n>0$ an that the same holds for any complex space whose 
universal covering is $\compl^N$ for some $N>0$. 

The first author wishes to thank Cales Casacuberta for having given
him the chance of visiting the Universitat de Barcelona and discussing
with him topics about localization of categories.

%%%%%%%%%%%%%%%%%%%%%%%%%%%%%%%%%%%%%%%%%%%%%%%%%%%%%%%%%%%%%%%%%%%%%%%%%
%%%%%%%%%%%%%%%%%%%%%%%%%%%%%%%%%%%%%%%%%%%%%%%%%%%%%%%%%%%%%%%%%%%%%%%%%%%
%%%%%%%%%%%%%%%%%%%%%%%%%%%%%%%%%%%%%%%%%%%%%%%%%%%%%%%%%%%%%%%%%%%%%%%%%%%

\section{Basic constructions}\label{base}

In this paper with $\mathbb{P}^n$ we will denote the $n$-th dimensional
projective complex space.
The general problem we are dealing with is to modify the category of complex
spaces to a category where the constant morphism ${\sf p}:\compl\to {\sf pt}$
 is invertibile. The task of inverting morphisms in a category, can be
accomplished by starting from an arbitrary category
 $\mathcal C$ with respect to a given family $S$ of morphisms
 satisfying suitable compatibility conditions
 (cfr. \cite{gz}). The category $S^{-1}\mathcal C$
 that we obtain is called the {\it localization} of $\mathcal C$ with
 respect to $S$. In this kind of generality, $S^{-1}\mathcal C$ is
not practical to work with.
 In this sense, reasonable categories are the ''homotopy categories''
associated to a
{\em model structure} in the sense of Quillen (i.e. endowed with a "good
 definition" of {\em weak equivalence} \cite{quillen}).

In this section we recall the main results of \cite{morvoe}. The
constructions made
there hold in the general context of a site with enough points in the sense of
 \cite{sga4}.  We restrict ourselves to the site $\mathcal S_T$
of complex spaces with the strong topology
%%%%%%%parte nuova
or that of schemes of finite type over a noetherian
 scheme $B$ of finite dimension, endowed with a Grothendieck topology which is
 weaker or as fine as the quasi compact flat topology.

%%%%%%%%%%%%%%%%%%%%%%%%%%%%%%%%%%%%%%%%%%%%%%%%%%%
%%%%%%%%%%%%%%%%%%%%%%%%%%%%%%%%%%%%%%%%%%%%%%%%%%%%

\subsection{Sheaves and simplicial objects: the categories $\fasc$ and
$\Delta^{\sf op}\fasc$}\label{colimiti}

Let $\mathcal S$ be the category of complex spaces or  schemes of
finite type over a noetherian scheme $B$. If we wish to do some kind
of homotopy theory on it, we should check the shape of colimits of
certain diagrams. Recall that given a diagram $\bf D$
\begin{equation}
\begin{split}
\xymatrix{A\ar[r]^i\ar[d]^f& X\\
B & }
\end{split}
\end{equation}
in a small category $\ccal$, an object $\mathrm{colim}_{\bf D}$ in $\ccal$ is the
colimit of $\bf D$ if and only if $\mathrm{colim}_{\bf D}$ fits in the
commutative diagram
\begin{equation}
\begin{split}
\xymatrix{A\ar[r]^i\ar[d]^f& X\ar[d]^p\\
B\ar[r]^g & \mathrm{colim}_{\bf D}}
\end{split}
\end{equation}
and ${\sf Hom}_{\ccal}(\mathrm{colim}_D, X)$ is the
limit of the diagram
 \begin{equation}
\begin{split}
\xymatrix{\homs_{\ccal}(A, Z)& \ar[l]_{i^*} \homs_{\ccal}(X,Z) \\
\homs_{\ccal}(B, Z)\ar[u]_{f^*} & }
\end{split}
\end{equation}
in the category of sets for any $Z\in\ccal$. In other words, this last
condition means
that ${\sf Hom}_{\ccal}(\mathrm{colim}_{\bf D}, X)$ are pairs of morphisms
$(\alpha, \beta)$, $\alpha\in\homs_{\ccal}(X,Z)$ and $\beta \in
\homs_{\ccal}(B,Z)$ with the property that $i^*\alpha=f^*\beta$.
The definition of colimit of an arbitrary diagram is similarly reduced
to the one of limit in the category of sets by applying
$\homs_{\ccal}(~~,Z)$. We are particularly interested in colimits of
diagrams of the kind
\begin{equation}
\begin{split}
\xymatrix{A\ar[r]^i\ar[d]^f& X\\
{\sf pt}& }
\end{split}
\end{equation}
where $i$ is an injection. In this paper, such colimits will sometimes
be called {\em quotient} of $X$ by
$A$ along $i$. In general, it may happen that the quotient does
not exist in the category $\mathcal{S}$ or if it exists, it is different from
the one taken in the underlying category of topological
spaces.
\begin{esempi}\label{esempio}
{\rm\begin{enumerate}
\item[1)] Let $\bf D$ be the diagrams
\begin{equation}
{\begin{split}
\xymatrix{\compl-0\ar[r]^i\ar[d]^f& \compl\\
{\sf pt}}&&\>\>\>\>
\xymatrix{\pro\ar[r]^i\ar[d]^f&\mathbb{P}^2\\
{\sf pt}&}
\end{split}}
\end{equation}
where $i$ are the canonical embeddings. Then the colimits of
${\bf D}$ in $\co$  are just a point in both cases, unlike their respective
colimits in the category of topological spaces.
\item[2)] Let $\bf D$ be the diagram
\begin{equation}
\begin{split}
\xymatrix{\mathbb Z\ar[r]^i\ar[d]^f&\compl\\
{\sf pt}& }
\end{split}
\end{equation}
where $i$ is the canonical injection. $\bf D$ has no colimit in $\mathcal
S$. Indeed, by contradiction, let $Z=\sf{colim}_{\bf D}$ in $\mathcal S$,
$p:\compl\to Z$
the corresponding canonical holomorphic function and $x=p(\mathbb
Z)$. Since there exists a non constant holomorphic function
$h:\compl\to\compl$ such that $h(n)=0$ for every
$n\in\mathbb Z$, $Z$ cannot be just the point $x$, moreover
$p^{-1}(x)=i(\mathbb Z)$. Let $U$ be a relatively compact
neighbourhood of $x$ and
$\{z^{(n)}\}\subset p^{-1}(U)$ a sequence with no accumulation
points. If $h:\compl\to\compl$ is a holomorphic function satisfying
$h(n)=0$ and $h(z^{(n)})=n$ for every $n\in\mathbb Z$, no holomorphic function
$g:Z\to\compl$ exists such that $g\circ p=h$.

A similar argument can be used to prove that the diagram
\begin{equation}
\begin{split}
\xymatrix{\compl\ar[r]^i\ar[d]^f&\compl\times\compl\\
{\sf pt}& }
\end{split}
\end{equation}
where $i$ is the injection $\compl\to\{0\}\times\compl$, has no
colimit in $\mathcal{S}$.
\item[3)] Let $\bf D$ be the diagram
\begin{equation}
\begin{split}
\xymatrix{\{0\}\cup\{1\}\ar[r]^i\ar[d]^f&\mathbb A^1_k\\
{\sf pt}& }
\end{split}
\end{equation}
where $\mathbb A^1_k$ is the affine line over a field $k$ and $i$ is
the embedding of the corresponding rational points. Then, since the
$k$-algebra of the polynomials
$P(x)$ of the form $a+x(x-1)Q(x)$, $a\in k$, is not finitely
generated, $\bf D$ has no colimit in the category of the algebraic schemes
of finite type over $k$.
\end{enumerate}}
\end{esempi}
We therefore enlarge $\mathcal{S}$ to a category which contains the
colimits of all diagrams and, at the same time, have a
``reasonably good'' shape from our point of view. Such a category is $\fasc$:
the objects are sheaves of sets on a site $\mathcal{S}$ endowed with
the Grothendieck topology $T$ and morphisms are maps of sheaves of sets.
Recall that a {\it sheaf of sets} on $\mathcal S_T$ (or an
arbitrary site) is a controvariant functor $\mathcal{F}:S_T\to Sets$
satisfying the following conditions:
\begin{enumerate}
\item $\fcal(\emptyset)=\{\sf{pt}\}$, where $\sf{pt}$ is the final
  object of $\mathcal{S}_T$;
\item let $q:E\to X$ be a covering for the topology $T$, $q_1$
and $q_2$ respectively the canonical projections $E\times_X E\to E$; then
\begin{equation}
\mathcal{F}(X)\stackrel{q^*}{\to}
\mathcal{F}(E)\stackrel{q_1^*}{\underset{q_2^*}{\rightrightarrows}}
\mathcal{F}(E\times_X E)
\end{equation}
is an exact sequence of sets i.e.
$$
q^*\mathcal{F}(X)=\{a\in\mathcal{F}(E):q_1^*(a)=q_2^*(a)\}
$$
\end{enumerate}

Let $\mathbf{Y}(X):={\sf Hom}_{\mathcal S}(\cdot,X)$. The functorial equality
$$
{\sf Hom}_{\mathcal S}(A,B)={\sf Hom}_{Funt(\mathcal S^{{\rm }op},
Sets)}(\mathbf{Y}(A),\mathbf{Y}(B))
$$
is known as Yoneda Lemma.
The Yoneda embedding is a faithfully full functor
$\mathbf{Y}:\mathcal{S}\hookrightarrow Fun(\mathcal{S}^{op}, Sets)$.
If the topology $T$ is not finer than the quasi compact flat topology,
then the image of $\mathbf{Y}$ is contained in the full subcategory
$\fasc$.

\begin{theorem}
Let $X\in S_T$ and $T$ a topology not finer than the quasi compact
flat topology or the strong topology in the holomorphic case. Then the
functor ${\sf Hom}_{\mathcal S}(\cdot, X)$ is a
sheaf for the topology $T$.
\end{theorem}
\demo In the algebraic case, we restrict the proof to the case in
which $\mathcal S_T$ is a site of shemes of finite type over a base as
we have been assuming from the very beginning. Then the
conclusion  follows from the theorem of Amitsur \cite{milne}.

Assume now that $\mathcal S_T$ is the site of complex spaces and let
$q:E\to Z$ be an open covering of the complex space $Z$ . Then, the sequence
\begin{equation}\label{seq}
E\times_Z E\stackrel{q_1}{\underset{q_2}{\rightrightarrows}}
E\stackrel{q}{\to} Z
\end{equation}
is exact as sequence of sets. We have to prove that the sequence of sets
\begin{equation}
{\sf Hom}_{\mathcal S}(E\times_Z
E,X)\stackrel{q_1^*}{\underset{q_2^*}{\leftleftarrows}}
{\sf Hom}_{\mathcal S}(E,X)\stackrel{q^*}{\gets}
{\sf Hom}_{\mathcal S}(Z,X)
\end{equation}
is exact, as well. Suppose that $q_1^*f=q_2^*f$ with $f\in
{\sf Hom}_{\mathcal S}(E,X)$. Since $q$ is continuous, surjective and
$Z$ has the quotient topology induced by $q$, applying the functor ${\sf
  Hom}_{Top}(\cdot, X)$ to the exact sequence (\ref{seq}) we obtain an
exact sequence, hence a
 continuous map $f':Z\to X$ such that $f=f'\circ q$. It follows that $f'$ is
 holomorphic, $f$ being holomorphic and $q$ a local biholomorphism.
\enddemo
\medbreak
\noindent
The category $\fasc$ is complete and cocomplete. Indeed, the limit
of a diagram $\bf D$ in
$\fasc$ is the functor $U\rightsquigarrow \lim{\bf D}(U)$ which is a
sheaf for the topology $T$. As for the colimit, it is defined as
$a_T(U\rightsquigarrow{\sf colim}{\bf D}(U))$ where $a_T$ is the
associated sheaf. In particular it possesses two canonical objects: an
initial sheaf $\emptyset$,
%%%%%%%%%%%%%%%parte nuova
the sheaf that associates the empty set to any
element of the site, except for the initial object of the site
$\mathcal{S}$  to which it
associates the one point set and the final sheaf, which we will
denote as $\pts$ or $Spec~B$ if the objects of the site are complex
spaces or schemes over $B$, respectively.

\medbreak\noindent

We now would like to consider the localized category $p^{-1}\fasc$,
where $p:\compl\to \sf{pt}$ (or $p:\aff\to Spec~k$). Moreover, we wish
the localized category to have supplementary structures such as
the ones we would get if $p^{-1}\fasc$ were equivalent to the homotopy
category of an appropriate model structure on $\fasc$.
Basically, a model structure on a category $\ccal$ is the data of three
classes of morphisms: weak equivalences, cofibrations and fibrations
satisfying five axioms CM1,$\cdots$, CM5 (see \cite{quillen}) with the
request that, in addition, the factorizations of CM5 are functorial.
We do not know about the existence of such model structure on
$\fasc$. This is a particular
case of the more general and complicated question on whether a
localized category $S^{-1}\ccal$ is equivalent to the homotopy
category associated to some model structure on $\ccal$.
Some results of this kind
are known in the case $\ccal$ itself is a homotopy category (see
\cite{bousfield}). To use them, we are forced to
embed $\fasc$ in the ``simplest'' category we know that is endowed of
a model structure, namely $\Delta^{op}\fasc$, the category of simplicial
objects in $\fasc$.

A {\it simplicial object } $\mathcal{X}$
%%%%%%%%%%%%parte nuova
in $\ccal$ is a sequence
$\{\mathcal{X}_i\}_{i\geq 0}$ of objects of $\ccal$ with a sequence
$\partial_i^n:\mathcal{X}_n\to\mathcal{X}_{n-1}$ of morphisms for $n\geq 1$,
$i=0,1,\cdots,n$ called {\em faces} and a sequence $\sigma_i^n:\mathcal{X}_n\to
\mathcal{X}_{n+1}$ of morphisms for $n\geq 0$, $i=0,1,\cdots,n$ called
{\em degenerations}, satisfying the following conditions
\begin{enumerate}
\item[1)] $\partial_i\partial_j=\partial_{j-1}\partial_i\hskip 30pt if~ i<j$
\item[2)] $\sigma_i\sigma_j=\sigma_{j+1}\sigma_i\hskip 30pt if~ i\leq j$
\item[3)] $\partial_i\sigma_j=\begin{cases}\sigma_{j-1}\partial_i & if~i<j\\
identity & if~ i=j~or~i=j+1\\
\sigma_j\partial_{i-1} & if~i>j+1
\end{cases}$
\end{enumerate}
A {\it morphism} $f:\xcal\to\ycal$ of two simplicial objects $\xcal=
\{\mathcal{X}_i\}_{i\geq 0}$, $\ycal=\{\mathcal{Y}_i\}_{i\geq 0}$ of
$\ccal$
is a sequence $\{f_i\}_{i\geq 0}$ of morphisms $f_i:\xcal_i\to\ycal_i$
which make the diagrams

$$
{\begin{split}
\xymatrix{\xcal_i\ar[r]^{f_i}\ar[d]^{\sigma^n_i}&\ycal_i
\ar[d]^{\sigma^n_i}\\
\xcal_{i+1}\ar[r]^{f_{i+1}}&\ycal_{i+1}}&&
\xymatrix{\xcal_i\ar[r]^{f_i}\ar[d]^{\partial^n_i}&\ycal_i
\ar[d]^{\partial^n_i}\\
\xcal_{i-1}\ar[r]^{f_{i-1}}&\ycal_{i-1}}
\end{split}}
$$
commutative.

\noindent With this notion of morphism, the family of simplicial
objects of $\ccal$ forms a
category denoted by $\Delta^{\sf op}\ccal$. Given $X\in\ccal$ we
denote by the
same symbol the {\it constant simplicial object} defined by
$\mathcal X_i=X$, $\partial_i^n=\sigma_i^n={\sf Id}_X$, for every $i,n$.

Suppose that $\ccal$ has a final object $*$, direct products and
direct coproducts. Let $[n]$ be the set $\{0,1,\cdots n\}$
Then,  for every integer $n\geq 0$,
denote by $\Delta[n]$ the simplicial object that at the level $m$ has
as many copies of $*$ as nondecreasing monotone functions $[m]\to [n]$.
The $m+1$ injective functions $[m-1]\to [m]$ induce the faces and the
$m$ surjective functions $[m]\to [m-1]$ induce the degeneracies of
$\Delta[n]$. On each copy of $*$ they act as the identity morphism.
Notice that in $\Delta[n]_n$ there is only one nondegenerate element,
namely the one corresponding to the identity. For example,
$\Delta[1]$ is described as $\Delta[1]_i=\amalg_{j=1}^{i+2}*$ for each
$i\geq 0$ and of the three $*$ in degree $1$, two of them are the
degenerations of of the $*$ in degree $0$. The two $*$ in degree zero
are the images through the face morphisms of the nondegenerate $*$ in
degree 1.

 \begin{oss}\label{star}
{\em For every simplicial object $\xcal$, the pair of $*$ in degree
  $0$ defines two morphisms $\epsilon_0$ and
  $\epsilon_1:\xcal\to\xcal\times\Delta[1]$.}
\end{oss}

Let $\xcal$, $\ycal$ two objects of $\Delta^{\sf op}\ccal$.
\begin{definition}\label{omot}
A {\it homotopy}
between two morphisms $f,g:\xcal\to\ycal$ is a morphism
$H:\xcal\times\Delta[1]
\to\ycal$ such that $H\circ\epsilon_0=f$, $H\circ\epsilon_1=g$.
\end{definition}
In particular, this definition gives a notion of homotopy for objects
and morphisms of $\ccal$.
\medbreak\noindent
\begin{esempi}\label{esempi25}
{\rm \begin{enumerate}
\item[1)] Let
$$
\Delta_{\sf top}^n=\{(t_0,t_1,\cdots,t_n)\in \mathbb{R}^{n+1}:0\leq t_i\leq 1,
\Sigma_i t_i=1\}.
$$
The collection $\{\Delta^n_{\sf top}\}_n$ forms a cosimplicial
topological space $\Delta^\bullet_{\sf top}$ with the standard coface
morphisms $\partial^i$ (inclusion of the face missing the vertex
$v_i$) and codegenerations $\sigma^i$ (proiection from $v_i$ on the
corrisponding face).
\item[2)] Let $\ccal$ be the category of sets. An object $A_\bullet=
\{A_i\}_{i\ge 0}$ of $\Delta^{\sf op}\ccal$ is called a {\it simplicial
set}. The {\em geometrical realization} of $A_\bullet$ is the topological
space
$$ |A_\bullet| =\dfrac{\amalg_n A_n\times\Delta^n_{\sf
top}}{(\partial_i(a),t)\sim (a,\partial^i(t)).} $$
A morphism
$\phi:A_\bullet\to B_\bullet$ of simplicial objects is said to be a {\em
weak equivalence} if its topological realization $|\phi|:|A_\bullet|\to
|B_\bullet|$ is a weak equivalence, i.e. the homorphisms
$|\phi|_*:\pi_k(|A|,a) \to\pi_k(|B|,|\phi|(a))$,  between the
homotopy groups are isomorphisms, for all $k> 0$ and a bijection for $k=0$.
\item[3)] Let ${\sf Top}$ be the category of topological spaces with
continuous maps. Then the functor ${\sf Sing}:{\sf Top}\to \Delta^{\sf
  op}\ins$, which associates to a
topological space $K$ the simplicial
%%%%%%parte nuova
set ${\sf Hom}_{\topo}
(\Delta^\bullet_{\sf top},K)$ is a functor that is left adjoint to
$A_\bullet\rightsquigarrow |A_\bullet|$. The pair of adjoint functors
$(\sf{Sing}, |~~|)$
  \begin{equation}
\xymatrix{\Delta^{\sf op}{\sf Ins}\ar@<.5ex>[r]^{\hskip
  15pt|~~|}& {\topo}\ar@<.5ex>[l]^{\hskip 15pt {\sf Sing}}}
\end{equation}
sends simplicial homotopies in the sense of Definition \ref{omot} to
homotopies of topological spaces and viceversa.
\end{enumerate}}
\end{esempi}
\noindent A simplicial object in $\fasc$ is said to be a {\it
simplicial sheaf}. For the time being, we will consider $\fasc$ as the
full subcategory of $\Delta^{\sf op}\fasc$, identified with constant
simplicial sheaves.

%%%%%%%%%%%%%%%%%%%%%%%%%%%%%%%%%%%%%%%%%%%%%%%%%%%%%%%%%%%%%%%%%%%%%

\subsection{Simplicial localization}

The following will endow $\Delta^{\sf op}\fasc$ with a model stucture
in the sense of Quillen:
\begin{definition}
 A morphism $f:\mathcal{G}\to\mathcal{F}$ of simplicial
sheaves is a {\it weak
equivalence} if for every point $x$ of a complex space or a scheme over $B$,
$f_x:\mathcal{G}_x\to \mathcal{F}_x$ is a weak equivalence of simplicial sets
($\mathcal{G}_x$ and $\mathcal{F}_x$ being the respective stalks over $x$ of
$\mathcal{F}$ and $\mathcal{G}$).
 \end{definition}

An injective  morphism $f:\xcal\to\ycal$ is
said to be a
%%%%%%parte nuova
{\it simplicial cofibration}.

 A {\em lifting} in a commutative square of morphisms

 \begin{equation}\label{diag}
\begin{split}
\xymatrix{\mathcal A\ar[r]^{q}\ar[d]^{j}& \xcal\ar[d]^{f}\\
\mathcal B\ar[r]^{r} & \ycal}
\end{split}
\end{equation}
is a morphism $h:\mathcal B\to \mathcal X$ which makes the diagram
commutative. In such situation we say that ${j}$ has the
{\em left lifting property} with respect to ${f}$ and
 ${f}$ has the {\em right lifting property} with respect to $j$.

A morphism $f:\mathcal X\to \mathcal Y$ is called a {\it fibration} if all
diagrams (\ref{diag}) admit a lifting, for all {\em acyclic
  cofibrations} $j$ (cofibration and weak equivalence simultaneously).

%%%%%%%%%%%%%%%%%%%%%%%%%%%%%%%%%%%%%%%%%%%%%%%%%%%%%%%%%%%%%%%%%%%%%%%%%%%

%%%%%%%%%%%%%%%%%%%%%%%%%%%%%%%%%%%%%%%%%%%%%%%%%%%%%%%%%%%%%%%%%%%%%%%%%

The classes of weak equivalences, cofibrations and fibrations give
$\Delta^{\sf{op}}\fasc$ a structure of {\em simplicial} model category
as shown in  \cite{jardine}. Under these assumptions, there exists a
localization of $\Delta^{\sf{op}}\fasc$ with respect of the weak
equivalences.
In other words, there exists a category which we will
denote by $\hsolo$ and a functor

$$
{\sf l}:\Delta^{\sf op}\fasc\to\hsolo
$$
which has the properties
\begin{itemize}
\item[1)]  if $f$ is a weak equivalence, ${\sf l}(f)$ is an
  isomorphism;
\medbreak
\item[2)] the property is universal, namely, if another category $\mathcal{C}$
exists and it is endowed with a functor ${\sf l'}:\Delta^{\sf op}\fasc
\to\mathcal{C}$ with the same
property as {\sf l}, then there exists a unique functor ${\sf u}:
\mathcal{H}_s^{olo}\to \mathcal{C}$ such that ${\sf l'}={\sf u}\circ {\sf l}$.
\end{itemize}

An object $\mathcal X$ of $\Delta^{\sf op}\fasc$
\begin{enumerate}
\item[1)] is called {\it cofibrant} if $\emptyset\to\mathcal{X}$ is a
cofibration;
\item[2)] is called {\it fibrant} if $\mathcal{X}\to {\sf pt}$ is a fibration.
\end{enumerate}

\subsection{Notations}\label{boh1}
 \begin{enumerate}
\item[1)] We denote $\sf{pt}$ the simplicial constant sheaf defined as
  the associated sheaf to the
the presheaf which associates to an object of $\mathcal S$ the set
consisting of one element.
The {\em pointed category} associated to $\Delta^{\sf op}\fasc$ is the
category $\Delta^{\sf op}_\bullet\fasc$ whose objects are the pairs
$(\mathcal{X},x)$ where
$\mathcal{X}\in\Delta^{\sf op}\fasc$ and $x:\sf{pt}\to X$ is a
morphism; a morphism
of pairs $(\mathcal{X},x)\to (\mathcal{Y},y)$ is a morphism $f:\xcal
\to\ycal$ such that $f\circ x =y$. As pointed sheaf, $\sf{pt}$ will
  stand for $(\sf{pt}, \sf{pt})$.

There is a pair of adjoint functors
\begin{equation}\label{nota}
\Delta^{\sf op}\fasc\underset{\sf t}{\stackrel{+}{\rightleftarrows}}
\Delta^{\sf op}_\bullet\fasc
\end{equation}
where $\sf t$ is the forgetful functor and ${\sf +}$ is defined by
$:\xcal
\rightsquigarrow\xcal_+$ with $\xcal_+:=\xcal\amalg{\sf pt}$, pointed
by $\sf{pt}$.
\item[2)] Let $f:\ycal\to\xcal$ be a morphism of (pointed) simplicial
sheaves. The symbol ${\sf cof}(f)$ denotes the colimit of the diagram
  $$
  \xymatrix{\ycal\ar[r]^f\ar[d] & \xcal\\{\sf pt }& }
  $$
 (pointed by the image of $\ycal$) where ${\sf pt}$ is a point.
${\sf cof}(f)$ is called the {\em cofibre} of $f$. If  $f$ is a cofibration
 the cofibre of $f$ is sometimes denoted by $\xcal/\ycal$.
\item[3)] Let $\xcal$ and $\ycal$ be pointed simplicial sheaves. The sheaf
$\xcal\vee\ycal$ is, by definition, the colimit of
  $$
  \xymatrix{{\sf pt}\ar[r]\ar[d] & \xcal\\
  \ycal & }
  $$
 pointed by the image of  ${\sf pt}.$
\item[4)]\label{smash} The pointed simplicial sheaf $\xcal\wedge\ycal$
is defined by  $\xcal\times\ycal/\xcal\vee\ycal$.
\item[5)] The simplicial pointed constant sheaf $S_s^1$ is defined
by $\Delta[1]/
\partial\Delta[1]$ where $\partial\Delta[1]$ is the simplicial subsheaf
  of $\Delta[1]$ costisting in the union of the images of the face
morphisms of $\Delta[1]$. For $n\in\mathbb N$ we set
$S^n_s=S^1_s\wedge\stackrel{n}{\cdots}\wedge S_s^1$.
\end{enumerate}

\medbreak\noindent

\begin{oss}
{\em Performing the same constructions as for $\Delta^{\sf op}\fasc$
  we obtain a homotopy category $\hsolop$.}
\end{oss}
\medbreak
For a more complete description of the main properties of $\hsolo$ and
$\hsolop$ we refer to \cite{quillen} e \cite{morvoe}. Here we only
recall
%%%%%%%%%parte nuova
a proposition that will be used later.

\begin{proposition}\label{mod-simpl}
Let $i:\mathcal{Y}\to\mathcal{X}$ be a simplicial cofibration of
pointed simplicial sheaves. Then, for every pointed simplicial sheaf
$\mathcal{Z}$, the morphism
$i$ induces a long exact sequence of pointed sets and groups
(see the proof of Lemma \ref{gruppo})
  \begin{multline}\label{sel}
{\sf Hom}_{\hsolop}(\mathcal{Y},\mathcal{Z})\stackrel{{i}^*}{\gets}
{\sf Hom}_{\hsolop}
(\mathcal{X},\mathcal{Z})\stackrel{\pi^*}{\gets} {\sf Hom}_{\hsolop}
(\mathcal{X}/\mathcal{Y},
\mathcal{Z})\gets\\ {\sf Hom}_{\hsolop}(\mathcal{Y}\wedge S^1_s,\mathcal{Z})
\stackrel{i^*}\gets {\sf Hom}_{\hsolop}(\mathcal{X}\wedge S_s^1)
\stackrel{\pi^*}{\gets} {\sf
  Hom}_{\hsolop}(\mathcal{X}/\mathcal{Y}\wedge S_s^1,\mathcal{Z}) \cdots
\end{multline}
\end{proposition}
%%%%%%%%%%%%%%%%parte nuova
This proposition is a particular case of Proposition $4'$ of \cite{quillen}.

The Yoneda embedding, induces a functor
%%%%%%%%parte nuova
$\yon_s:\mathcal{S}_T\to\hsolo$ which
is a full embedding (see the Proposition 1.13 and Remark 1.14 of
\cite{morvoe}). However, in general, it is more difficult to describe the
morphisms betweeen objects in $\hsolo$. Indeed, ${\sf Hom}_{\hsolo}
(\mathcal{Y},\mathcal{X})$ is obtained as a quotient of the set of diagrams
$$
\mathcal{Y}\underset{\gets}{\stackrel{\sf s}{\sim}}\mathcal{Y}'\to
  \mathcal{X}
$$
of $\Delta^{\sf op}\fasc$ where $s$ is a weak equivalence.

$\hsolo$ (or its pointed version) is the appropriate category in which
we are going to invert $p:\compl\to{\sf pt}$. In the next section we
will give a model structure to $\Delta^{\sf op}\fasc$ whose weak
equivalences contain $p$, and are in a sense the ``smallest'' class
containing all the base changements of $p$ as well. Such weak equivalences
are written in terms of morphisms in $\hsolo$ and the homotopy category
associated to this model structure is the localization of $\hsolo$
with respect to the weak equivalences.

%%%%%%%%%%%%%%%%%%%%%%%%%%%%%%%%%%%%%%%

\subsection{Affine localization}\label{affine}

%%%%%%%%%%parte nuova
Unless otherwise mentioned, the results presented in this subsection are
taken from section 3.2 of \cite{morvoe}.

\begin{definition}\label{a1loc}
A simplicial sheaf $\mathcal{X}\in\Delta^{\sf op}\fasc$ is said to be {\it
$\aff$-local} (or $\compl$ local in the complex case) if the projection
$\mathcal{Y}\times\aff\to\mathcal{Y}$ induces a bijection of sets
$$
{\sf Hom}_{\hsolo}(\mathcal{Y},\mathcal{X})\to
{\sf Hom}_{\hsolo}(\mathcal{Y}\times\aff,\mathcal{X})
$$
\noindent for every $\mathcal{Y}\in\Delta^{\sf op}\fasc$.
\end{definition}
In what follows we describe a new structure of models on $\Delta^{\sf
  op}\fasc$, which we will call {\em affine}.

A morphism $f:\mathcal{X}\to\mathcal{Y}$ is called:
\begin{enumerate}
\item [1)] an {\em affine (or $\aff$ in the algebraic case or $\compl$
  in the complex case) weak equivalence} if, for every $\aff$-local
simplicial sheaf $\mathcal{Z}\in\Delta^{\sf op}\fasc$
$$
f^*:{\sf Hom}_{\hsolo}(\mathcal{X},\mathcal{Z})\to
{\sf Hom}_{\hsolo}(\mathcal{Y},\mathcal{Z})
$$
is a bijection;
\item[2)] an {\it affine cofibration} if it is  injective;
\item[3)] an {\em affine fibration} if all diagrams (\ref{diag})
admit a lifting,  where $j$ is any affine cofibration and affine
weak equivalence.
\end{enumerate}
An object $\xcal$ of $\Delta^{\sf op}\fasc$ is called
\begin{enumerate}
\item[1)] $\aff$-{\it fibrant} if the canonical morphism $X\to*$ is an affine
fibration;
\item[2)] $\aff$-{\it cofibrant} if $\emptyset\to{X}$ an affine
 cofibration.
\end{enumerate}

%%%%%%%%%%%%%%%%%%%%%%%%parte nuova
\begin{theorem}(cfr. Theorem 3.2, \cite{morvoe})
The structures listed above endow $\Delta^{op}\fasc$ of a model
structure, which will be called {\em affine model structure} or {\em
  $\aff$ model structure}.
\end{theorem}

The localized category with respect of the affine weak equivalences is denoted
as $\holo$ and its pointed version as $\holop$.

\begin{oss}\label{pippo}
{\em \begin{enumerate}
\item Any object of $\Delta^{\sf op}\fasc$ is both (simplicially)
  cofibrant and  $\aff$-cofibrant.
\item\label{vai} If $f:\ycal\to\xcal$ is a simplicial weak equivalence
(respectively  a simplicial cofibration) then
  it is an affine weak equivalence (respectively an affine
  cofibration). Therefore, the affine localization functor
  $\Delta^{\sf op}\fasc\to \holo$ factors as $\Delta^{\sf
    op}\fasc\to\hsolo\to \holo$, where the first functor is the
  simplicial localization and the second is the identity on objects
  and identity on the fractions representing morphisms. However, the
  functor $\hsolo\to \holo$ is {\em not} an equivalence of categories.
\item The same classes of pointed morphisms, give $\Delta^{\sf
  op}_\bullet\fasc$ a model structure.
\end{enumerate}}
\end{oss}

%%%%%%%%%%%%%%%parte nuova

 Proposition \ref{mod-simpl} holds for $\holo$ as well
\begin{proposition}\label{mod-aff}
Let ${j}:\mathcal{Y}\to\mathcal{X}$ be an affine cofibration (i.e an
injection of simplicial pointed sheaves). Then, for every simplicial
pointed sheaf $\mathcal{Z}$, the morphism $j$ induces long exact
sequence of pointed sets and groups
\begin{multline}
{\sf Hom}_{\holop}(\mathcal{Y},\mathcal{Z})\stackrel{{j}^*}{\gets}
{\sf Hom}_{\holop}(\mathcal{X},
\mathcal{Z})\stackrel{\pi^*}{\gets} {\sf Hom}_{\holop}(\mathcal{X}/\mathcal{Y},
\mathcal{Z})\gets\\{\sf Hom}_{\holop}(\mathcal{Y}\wedge S^1_s,\mathcal{Z})
\stackrel{{j}^*}{\gets}{\sf Hom}_{\holop}(\mathcal{X}\wedge S_s^1)
\stackrel{\pi^*}{\gets}
{\sf Hom}_{\holop}(\mathcal{X}/\mathcal{Y}\wedge S_s^1,\mathcal{Z}) \cdots
\end{multline}
\end{proposition}
%%%%%%%%%%%%parte nuova
The proof of such a statement is the same as for the Proposition
\ref{mod-simpl}.

\subsection{Hyperbolic simplicial sheaves}
Let us go back to the concept of hyperbolicity.
\begin{definition}\label{spazio-iper} A simplicial sheaf $\xcal$ is said to be
{\em hyperbolic} if it is $\aff$-local. Let $\ccal$ be a simplicial subsheaf of
$\xcal$. The simplicial sheaf $\xcal$ is said to be
{\em hyperbolic mod $\ccal$} if $\xcal/\ccal$ is hyperbolic.
\end{definition}

\begin{definition}\label{ris-iper}
A {\it hyperbolic resolution} of $\xcal$ is a morphism of simplicial sheaves
$\mathfrak{r}:\xcal\to {\widetilde\xcal}$ where $\widetilde{\xcal}$ is a
hyperbolic simplicial sheaf and $\mathfrak r$ is an affine weak equivalence.
\end{definition}

A {\it hyperbolic resolution functor} is a pair $(\mathfrak{I}, \mathfrak{r})$
where $\mathfrak{I}$ is a functor
  $$
  \Delta^{\sf op} \fasc\to\Delta^{\sf op}\fasc
  $$
and $\mathfrak{r}$ is a natural transformation ${\sf Id}\to\mathfrak{I}$ such that every morphism $\xcal\to\mathfrak{I}
(\xcal)$ is a hyperbolic resolution.

From Proposition 2.19 of \cite{morvoe} we derive the following, fundamental
result:
\begin{theorem}\label{localizz}
There exists a hyperbolic resolution functor
$(\mathfrak{Ip},\mathfrak{r})$ with the following properties:

\begin{enumerate}
\item[1)] for every $\mathcal X\in\Delta^{\sf op}\fasc$ the simplicial
  sheaf $\mathfrak{Ip}(\xcal)$ is hyperbolic and (simplicially)
  fibrant;
\item[2)] $\mathfrak{r}$ is an affine equivalence and a cofibration;
\item[3)]\label{imm} let $\mathcal{H}_{s,\aff}$ be the full
subcategory in $\hsolo$ of $\aff$-local (hyperbolic) objects.
$\mathfrak{Ip}$ sends an affine weak equivalence to a
simplicial weak equivalence, hence it induces a functor
$L:\hsolo\to \mathcal{H}_{s,\aff}$, that factors as $\hsolo\to \holo\to
  \holo_{s,\aff}$, where the first functor is the identity on objects
  (see also Remark \ref{pippo} (\ref{vai}));
\item[4)] the canonical immersion
  $I:\mathcal{H}_{s,\aff}\hookrightarrow \hsolo$
is a right adjoint of $L$ .
\end{enumerate}
Furthermore, $\mathcal{H}_{s,\aff}$ is a category equivalent to
$\holo$.
\end{theorem}
\medbreak
Given $X=\xcal\in \boldsymbol{\mathcal{F}}_T$, $\mathfrak{Ip}(X)$
is the hyperbolic simplicial sheaf associated to the simplicially
constant sheaf $X$. However, due to its rather
involved construction, the use of $\mathfrak{Ip}(\xcal)$ is
problematic even in the case when $X$ is a complex space or a scheme
over $k$.

Therefore, in general, the previous result shall be considered as an
existence theorem. Nevertheless, it may occur that, in some particular
cases, the class of $\mathfrak{Ip}(\xcal)$ in $\holo$ could be
represented by an understandable object, or even by a hyperbolic
space (e.g. \ref{esem}). In order to give a more
precise idea of the difficulties involves, let $\underline{{\sf
Hom}}(\xcal,\zcal)$ be the right adjoint functor to $\ycal
\rightsquigarrow\ycal\times\xcal$, where the objects are simplicial
sheaves. Let us define ${\sf Sing}_\bullet^{\aff}(\xcal)$ to be the simplicial
sheaf $\{\underline{{\sf Hom}}(\Delta^n_{\aff}, \xcal_n)\}_{n\ge 0}$, where
$\Delta_{\aff}^\bullet$ is the cosimplicial sheaf such that $\Delta^n_{\aff}=
\mathbb{A}^n$ for every $n$ and the structure morphisms are as
described in page
88 of \cite{morvoe}. Then, the class $\mathfrak{Ip}(\xcal)$ is
defined to be the simplicial sheaf $$ ({\sf Ex}\circ{\sf
Sing}_\bullet^{\aff})^{\mathbf\omega}\circ {\sf Ex}(\xcal) $$ where
$\xcal\to {\sf Ex}(\xcal)$ is a fibrant simplicial resolution and
$\mathbf{\omega}$ is a sufficiently large ordinal.

 \medbreak

 We conclude
this section by a short discussion on morphisms in localized categories.
Morphisms in a localized  category  ${S}^{-1}\ccal$ can be expressed
in terms of morphisms of $\ccal$ using the so called {\it calculus
of fractions}. More precisely, $$ {\sf Hom}_{ S^{-1}\ccal}(\mathcal X,
\mathcal Y)=\dfrac{\{\mathcal X\stackrel{s}{\gets}\mathcal X'\stackrel{f}{\to}
\mathcal Y: s\in S, \; f\in {\sf Hom}_{\ccal}(\mathcal X',\mathcal
Y)\}}{\sim} $$ where the elements of the numerator are called {\it
fractions} and $\sim$ is an equivalence between fractions.

If the localization is associated to a model structure $\ccal$
(as it happens for $\hsolo$ and $\holo$),we know that there are objects
$\mathcal X$, $\mathcal Y$ such that, ${\sf
  Hom}_{{S}^{-1}\ccal}(\mathcal X,\mathcal Y)$ is a quotient of
${\sf Hom}_{\ccal}(\mathcal X,\mathcal Y)$; for instance, if $\mathcal X$ is
cofibrant and $\mathcal Y$ is fibrant. Under these assumptions it can be proved
 that
\begin{equation}\label{quoz}
{\sf Hom}_{S^{-1}\ccal}(\mathcal X,\mathcal Y)=\dfrac{{\sf
    Hom}_{\ccal}(\mathcal X,\mathcal Y)}{\sim_{l}}
\end{equation}
where $f\sim_{l} g$ if and only if the morphism $f\amalg g$ factors
through a {\em cylinder object} ${\sf Cyl}(\xcal)$:
\begin{equation}
\begin{split}
\xymatrix{{\sf Cyl}(\xcal)\ar[dr]^{{\sf can}} & \\ \xcal\amalg \xcal\ar[u]
\ar[r]^{\hskip 8pt f\amalg g} & \ycal}.
\end{split}
\end{equation}
We recall that a cylinder associated to an object $X$ of a category $\ccal$
endowed with a model structure is an object ${\sf Cyl}(X)\in \ccal$
with morphisms $$ X\amalg X\stackrel{\sf i}{\to}{\sf Cyl}(X)\stackrel{\sf
can}{\to} X $$ such that ${\sf can}\circ{\sf i}={\sf can}\circ{\sf
id}_X\amalg {\sf id}_X$ and ${\sf can}$ is a weak equivalence. Cylinder
objects always exists in a category $\ccal$ endowed with a model
structure. Furthermore the morphism ${\sf i}$ can be chosen to be a
cofibration, ${\sf can}$ a fibration and the corrispondence $X\rightsquigarrow
{\sf Cyl}(X)$ functorial.

If $\ccal=\Delta^{\sf op}\fasc$,  a cylinder object for the affine
model structure  associated to
$\xcal$ may be taken to be $\mathcal X\times\aff$ where $\sf i$ is the
morphism $X\amalg X\to \mathcal X \times\aff$ determinated by the
inclusions at $0$ and $1$ (i.e. by the morphisms
$\mathcal X\to\mathcal X\times \{0\}$, $\mathcal X\to\mathcal X\times
\{1\}$) and
${\sf can}$ the projection onto $\mathcal X$.
We already observed that every
object of $\Delta^{\sf op}\fasc$ is cofibrant for both the model structures
on $\Delta^{\sf op}\fasc$. Consequently, if $\mathcal Y$ is fibrant
(respectively simplicially fibrant),
${\sf Hom}_{\holo}(\mathcal X, \mathcal Y)$ (respectively
${\sf Hom}_{\hsolo}(\mathcal X,\mathcal Y)$) is a quotient set of
${\sf Hom}_{\Delta^{\sf op}\fasc}(\mathcal X,
\mathcal Y)$. In the sequel, this fact will be extensively used.

%%%%%%%%%%%%%%%%%%%%parte nuova

\begin{lemma}
For any simplicial sheaf $\xcal$, the morphism $\mathfrak{r}:\xcal\to
\mathfrak{Ip}(\xcal)$ is universal in the category $\holo$
(respectively in the category $\hsolo$) in the
following sense: for any hyperbolic object $\ycal$
and morphism $f:\xcal\to \ycal$ in $\holo$,
there exists a unique
morphism $\tilde{f}:\mathfrak{Ip}(\ycal)\to \ycal$ in $\hsolo$.
\end{lemma}

\demo
Consider the commutative square
\begin{equation}
\begin{split}
\xymatrix{\xcal\ar[r]^{\mathfrak{r}_\xcal}\ar[d]^f &
  \mathfrak{Ip}(\xcal)\ar[d]^{\mathfrak{Ip}(f)}\\
\ycal\ar[r]^{\mathfrak{r}_\ycal}_{\cong} & \mathfrak{Ip}(\ycal).}
\end{split}
\end{equation}
By definition of $\aff$ weak equivalence, $\mathfrak{r}_\ycal$ is a
simplicial weak equivalence, since both $\ycal$ and
$\mathfrak{Ip}(\ycal)$ are hyperbolic (i.e. $\aff$
local). The map $\tilde{f}$ is defined as $\mathfrak{r}_\ycal^{-1}\circ
\mathfrak{Ip}(f)$. Note that $\mathfrak{Ip}(f)$ is a morphism in
$\hsolo$.
\enddemo

\begin{corollary}
Let $X$ and $Y$ be shaves with $Y$ hyperbolic and $f:X\to Y$ be a
morphism of sheaves. Then the composition $\tilde{f}\circ
\mathfrak{r}$ is a morphism of sheaves and the commutativity of the diagram
\begin{equation}
\begin{split}
\xymatrix{X\ar[r]^{\mathfrak{r}}\ar[d]^f &
  \mathfrak{Ip}(X)\ar[dl]^{\tilde{f}} \\
Y & }
\end{split}
\end{equation}
is in the category of sheaves, i.e. it is strictly commutative and not
only ``up to homotopy'' in $\hsolo$.
\end{corollary}

\demo
By the previous lemma, we have commutativity in the category $\hsolo$.
Remark 1.14 of \cite{morvoe} implies that $\homs(X,
Y)=\homs_{\hsolo}(X,Y)$ since both $X$ and $Y$ have simplicial
dimension zero. Therefore, equality of the morphisms $f$ and
 $\tilde{f}\circ\mathfrak{r}$ in $\hsolo$ is an equality of morphisms
of sheaves.
\enddemo

%%%%%%%%%%%%%%%%%%%%%%%%%%%%%%%%%%%%%%%%%%%%%%%%%%%%%%%%%%%%%%%%%%%%%%%
%%%%%%%%%%%%%%%%%%%%%%%%%%%%%%%%%%%%%%%%%%%%%%%%%%%%%%%%%%%%%%%%%%%%%%%

\section{Hyperbolicity and Brody hyperbolicity}
In this section we will compare the different notions of hyperbolicity
that we have introduced above. In
particular, we prove that a simplicial sheaf represented by a complex
space $X$ is hyperbolic if and only if $X$ is Brody hyperbolic.
This is a corollary of the following
\begin{theorem}\label{fibering}
A sheaf $X\in\fasc$ is a hyperbolic sheaf if and only if the
projection $U\times\aff\to U$ induces a bijection
$$
{\sf Hom}_{\fasc}(U,X)\to {\sf Hom}_{\fasc}(U\times\aff,
X)
$$
for every object $U\in\mathcal S_T$. Moreover, under this hypothesis, for
 every $Y\in\fasc$ there exists a bijection
\begin{equation}\label{vaff1}
{\sf Hom}_{\holo}(Y,X)\cong {\sf Hom}_{\hsolo}(Y,X)\cong {\sf
  Hom}_{\fasc}(Y,X).
\end{equation}
\end{theorem}
\begin{oss}
{\em If in (\ref{vaff1}) the sets have a group structure induced (up to
homotopy) by a group structure on $Y$ or by a cogroup structure (up to
homotopy) on $X$, the bijection is a group isomorphism.}
\end{oss}
Before beginning the proof, we fix, by the following commutative diagram
\begin{equation}\label{diagr}
\begin{split}
\xymatrix{&& & \hsolo\ar@<.5ex>[dd]^L\\
S_T\ar[urrr]^{\yon_s}\ar[drrr]_{\yon_{\aff}}
\ar[r]& \fasc\ar[r]^{\hskip -10pt {\sf cost}}&
  \Delta^{\sf op}\fasc\ar[ur]^{L_s}\ar[dr]_{L_{\aff}}\\
&&& \hskip 20pt\holo_{s,\aff}\cong\holo\ar@<.5ex>[uu]^I}
\end{split}
\end{equation}
the names of the functors involved in the proof. Notice that,
the first functor on the
left is the Yoneda embedding and $L$ are the localization functors.\medbreak
\noindent
{\bf Proof of Theorem \ref{fibering}.}
First of all we have the following bijections of sets
\begin{equation}
{\sf Hom}_{\holo}(\ycal, \xcal)\cong {\sf Hom}_{\holo}(L(\ycal),\xcal)\cong
{\sf Hom}_{\holo_{s,\aff}}(L(\ycal),\xcal)\cong {\sf Hom}_{\hsolo}(\ycal,\xcal).
\end{equation}
The left end side bijection is a consequence of the fact that the
canonical morphism
$\mathcal Y\to L(\mathcal Y)$ is an affine equivalence, the second one
follows from the
equivalence between $\holo$ and $\holo_{s,\aff}$ (Theorem
\ref{localizz}) by definition of hyperbolicity of a simplicial sheaf.
Finally, the
third one follows from the fact that $(L,I)$ is a pair of adjoint
functors (Theorem \ref{localizz}.(4)).

Assume now that $\xcal=X$ and $\ycal=Y$ are sheaves. Using the results
quoted in
\cite[Remark 1.14, p. 52 ]{morvoe} one shows that
\begin{lemma}\label{uguali}
Let  $X$, $Y$ be sheaves. Then
$$
{\sf Hom}_{\hsolo}(Y,X)\cong
{\sf Hom}_{\fasc}(Y,X).
$$
 \end{lemma}
\noindent This result implies the second assertion of Theorem
\ref{fibering}. Indeed, if $X$ is a hyperbolic sheaf, from Lemma
\ref{uguali} combined with the above considerations we get
$$
{\sf Hom}_{\holo}(Y,X)={\sf Hom}_{\hsolo}(Y,X)={\sf Hom}_{\fasc}(Y,X)
$$
for every sheaf $Y$. Moreover, Lemma \ref{uguali} also implies
the first assertion in the following weaker form: given a sheaf $X$,
the projection $U\times\aff\to U$ induces a bijection
$$
{\sf Hom}_{\hsolo}(U,X)\to {\sf Hom}_{\hsolo}(U\times\aff,
X)$$
for every $U\in S_T$ if and only if it induces a bijection
$$
{\sf Hom}_{\fasc}(U,X)\to {\sf Hom}_{\fasc}(U\times\aff,X).
$$
Thus, in order to finish the proof of Theorem \ref{fibering} in the
general case,
it is sufficient to prove the following: for every $U\in S_T$, the projection
$U\times\aff\to U$ induces a bijection
 $$
{\sf Hom}_{\hsolo}(U,X)\to
{\sf Hom}_{\hsolo}(U\times\aff, X)
$$
if and only if $X$ is hyperbolic. For this we use the following  Lemma 1.16 of
\cite{morvoe}:
\begin{lemma}\label{palle}
Let $\Sigma$ a set of objects of $\fasc$ such that, for every $U\in
\mathcal S_T$,
there exists an epimorphism $F\to U$ where $F$ is a direct sum of elements of
$\Sigma$. Then there exist a functor
$$
\Phi_{\Sigma}:\Delta^{\sf op}\fasc \to\Delta^{\sf op}\fasc
$$
and a natural transformation $\Phi_{\Sigma}\to {\sf Id}$ with the following
properties: given $\ycal$
\begin{enumerate}
\item[1)] for every $n\geq 0$ the sheaf of sets $\Phi_{\Sigma}(\ycal)_n$ is a
direct sum of sheaves belonging to $\Sigma$;
\item[2)] the morphism $\Phi_{\Sigma}(\ycal)\to \ycal$ is both a
(simplicial) weak
equivalence and a local fibration (i.e. the morphism induces on the
stalks a Kan fibration of simplicial sets).
\end{enumerate}
\end{lemma}
Since we will refer often to this lemma, we are going to recall here
how to construct the functor $\Phi_\Sigma$. Let $f:\xcal\to\ycal$ be a
morphism of simplicial sheaves; define $\Psi_{\Sigma,f}$ as the
colimit of
\begin{equation}
\begin{split}
\xymatrix{\amalg_{{\bf D}_n,n\geq 0} F\times \partial\Delta[n]\ar[r]\ar[d] & \xcal\\
\amalg_{{\bf D}_n,n\geq 0} F\times \Delta[n] &}
\end{split}
\end{equation}
where ${\bf D}_n$ is the set of the commutative squares of the kind
\begin{equation}
\begin{split}
\xymatrix{F\times\partial\Delta[n]\ar[r]\ar[d] & \xcal\ar[d]^f\\
F\times\Delta[n]\ar[r] &\ycal}
\end{split}
\end{equation}
and $F\in \Sigma$. Let $\alpha_1:\Psi_{\Sigma,f}\to \ycal$ be the
canonical morphism and $\Phi^{m+1}_{\Sigma,f}$ be
$\Psi_{\Sigma,\alpha_m}$. The $\Phi^i_{\Sigma,f}$ form a direct
system of cofibrations $\{\Phi^1_{\Sigma,f}\cdots\hookrightarrow
\Phi^i_{\Sigma,f}\hookrightarrow \Phi^{i+1}_{\Sigma,f}\hookrightarrow
\cdots\}$ whose colimit we will denote as $\Phi_{\Sigma,f}$. Such a
simplicial sheaf factors functorially $f: \xcal\to\Phi_{\Sigma,f}\to
\ycal$. The functor that associates to a simplicial sheaf $\ycal$ the
simplicial sheaf $\Phi_{\Sigma,\emptyset\to\ycal}$ satisfies the
properties of the Lemma.

\medbreak

In view of the Yoneda Lemma, we see that we can take as $\Sigma$ the
sheaves represented by objects in $S_T$. Indeed, for every sheaf
$W$, we have a surjective morphism (even as presheaves)
\begin{equation}
\amalg_{U\in S_T}\amalg_{s\in W(U)}{\sf Hom}_{\fasc}(\cdot, U)\to W.
\end{equation}
Thus, by Lemma \ref{palle},  given an arbitrary simplicial sheaf $\ycal$ there
exists $\ycal'$ such that $\ycal=\ycal'$ in
$\hsolo$ and $\ycal'_n=\amalg_{n_i}U_{n_i}$ with $U_{n_i}\in S_T$.
By definition, if $X$ is hyperbolic (i.e. $\aff$-local),
the projection ${\sf p}:\ycal\times \aff\to\ycal$ induces a bijection
$$
{\sf p}^*:{\sf Hom}_{\hsolo}(\ycal, X)\to
{\sf Hom}_{\hsolo} (\ycal\times\aff, X)
$$
for every $\ycal\in \fasc$. In particular, this holds if
$\ycal=U\in\mathcal S_T$, so using
Lemma \ref{uguali} we conclude that for every $U\in \fasc$
\begin{equation}\label{bij}
{\sf p}^*:{\sf Hom}_{\fasc}(U,X)\to{\sf Hom}_{\fasc}(U\times \aff, X)
\end{equation}
is a bijection.

Conversely, assume that (\ref{bij}) is a bijection for every $U\in
S_T$. Then, for every $U\in S_T$,
$$
{\sf p}^*:{\sf Hom}_{\hsolo}(U,X_+)\to {\sf Hom}_{\hsolo}(U\times\aff, X_+)
$$
is a bijection. Let $sk_n\ycal$ be the simplicial sheaf defined by
$$
sk_n\ycal=\begin{cases}(sk_n\ycal)_i=\ycal_i,& {\rm se}\ i\le n,\\
(sk_n\ycal)_i=\amalg_{\sigma_{u_j}\;\mathrm{degenerations}}
 \sigma_{u_1}\circ\cdots\circ\sigma_{u_{i-n}}(\ycal_n), &{\rm se}\
 i>n.
\end{cases}
$$
Such an object is called the {\em n-skeleton} of $\ycal$. The immersion
$$
i_n:sk_{n-1}\ycal\hookrightarrow sk_n\ycal
$$
is a cofibration and for $\ycal=\ycal'$, the cofibre $sk_n\ycal'/
sk_{n-1}\ycal'$ is isomorphic to the sheaf $\amalg_{n_i}U_{n_i}\wedge S^n_s$.

We use the following cofibration sequences:
\begin{equation}\label{uno}
sk_{n-1}\ycal'_+\to sk_n\ycal'_+\to \vee_{n_i} {U_{n_i}}_+
\wedge S_s^n\to  sk_{n-1}\ycal'_+\wedge S_s^1
\end{equation}
\begin{multline}\label{due}
\vee_n sk_n\ycal'_+\to \mathrm{\sf dirlim}_n sk_n\ycal'_+=\ycal'_+\to \\\vee_n
  (sk_n\ycal'_+\wedge S_s^1)\stackrel{{\sf id}-\vee i_n\wedge {\sf id}_{S^1}}
{\longrightarrow}\vee_n(sk_n\ycal'_+\wedge S_s^1).
\end{multline}
Notice that we are forced to take separate base points, since in the
algebraic case, we cannot assume that a simplicial sheaf $\zcal$ can be
considered as a pointed simplicial sheaf.

If $n=1$ the sequence (\ref{uno}) becomes
\begin{equation}
(\amalg_{0_i}U_{0_i})_+\hookrightarrow
 sk_1\ycal'_+\to \vee_{1_i}{U_{1_i}}_+\wedge S^1_s\to
 (\amalg_{0_i}U_{0_i})_+\wedge S_s^1.
\end{equation}
Thus the following  sequence
\begin{equation}
(\amalg_{0_i}U_{0_i}\times\aff)_+\hookrightarrow
 (sk_1\ycal'\times\aff)_+\to (\vee_{1_i}{U_{1_i}}
\wedge S^1_s)\wedge\aff_+\to
 \amalg_{0_i}{U_{0_i}}_+\wedge S_s^1\wedge\aff_+
\end{equation}
is a cofibration sequence as well. The projection ${\sf p}:\aff\to {\sf pt}$
maps the latter sequence to the former.

Applying ${\sf Hom}_{\hsolop}(~~, X_+)$ we get the long exact sequence
of pointed sets and groups

\begin{multline}\label{sel1}
 {\sf Hom}_{\hsolop}((\amalg_{0_i}{U_{0_i}})_+,X_+)\gets
 {\sf Hom}_{\hsolop}(sk_1\ycal'_+,X_+)\gets\\
 {\sf Hom}_{\hsolop}(\vee_{1_i}{U_{1_i}}_+\wedge S_s^1,X_+)
\gets{\sf Hom}_{\hsolop}((\amalg_{0_i}U_{0_i})_+\wedge S_s^1,X_+) \cdots
\end{multline}
as a particular case of the exact sequence (\ref{sel}). The morphism
${\sf p}^*$ induces maps from the sequence (\ref{sel1}) to the one
corresponding to $\ycal\times\aff_+$. We are
going to prove that ${\sf p}^*$ is a bijection of pointed sets, from
$$
{\sf A}={\sf Hom}_{\hsolop}((\amalg_{0_i}{U_{0_i}})_+,X_+),
$$
$$
{\sf B}={\sf Hom}_{\hsolop}(\vee_{1_i}{U_{1_i}}_+\wedge S_s^1,X_+),
$$
$$
{\sf C}={\sf Hom}_{\hsolop}(\amalg_{0_i}{U_{0_i}}_+\wedge S_s^1,X_+).
$$
\noindent
${\sf p}^*$ is bijective from $\sf A$, because by the adjunction
(\ref{nota}), we get
$$
{\sf A}={\sf Hom}_{\hsolo}
(\amalg_{0_i}U_{0_i}, X_+).
$$
Since direct sums of classes in $\hsolo$ are represented by direct
sums in $\Delta^{\sf op}\fasc$, we have that ${\sf A}=\amalg_{0_i} {\sf
  Hom}_{\hsolo}(U_{0_i}, X_+)$ and we conclude by using the assumption
we have on $X$.

\noindent Regarding the pointed set $\sf B$ we argue as follows: a
fibrant model of $X_+$
is of the kind $\tilde{X}_+$, where $\tilde{X}$ is a fibrant model of
$X$, thus it is a nonconnected simplicial sheaf. On the other hand,
$\vee_{1_i}{U_{1_i}}_+\wedge S_s^1$ is a pointed connected simplicial
sheaf. Since ${\sf B}$ is a quotient set of ${\sf Hom}_{\Delta^{\sf
    op}_\bullet\fasc}(\vee_{1_i}{U_{1_i}}_+\wedge S_s^1,
\tilde{X}_+)$, we conclude that ${\sf B}=*$, the constant map to the
base point,  because ${\sf Hom}_{\Delta^{\sf
    op}_\bullet\fasc}(\vee_{1_i}{U_{1_i}}_+\wedge S_s^1,
\tilde{X}_+)$ is. The same argument works for

$${\sf Hom}_{\Delta^{\sf op}_\bullet\fasc}(\vee_{1_i}{U_{1_i}}_+\wedge
S_s^1\wedge (\aff_+),\tilde{X}_+).$$
Thus $\sf p^*$ is an isomorphism on $\sf B$. The same argument shows
that $\sf p^*$ is also an isomorphism on $\sf C$. By the Five Lemma, we conclude that $\sf p^*$ is an
isomorphism from  ${\sf Hom}_{\hsolop}(sk_1\ycal'_+, X_+)$.

Similarly, we prove that the cofibration exact sequences  (\ref{uno})
yield that ${\sf  p}^*$ is bijective from ${\sf Hom}_{\hsolop}
(sk_n\ycal'_+,X_+)$, for every $n\geq 0$.

Since ${\sf p}^*$ is bijective from ${\sf C}={\sf Hom}_{\hsolop}
((\amalg_{0_i}U_{0_i})_+\wedge S_s^1, X_+)$ the same holds for
${\sf Hom}_{\hsolop}(sk_1\ycal'_+\wedge S_s^1,X_+)$ and consequently for
${\sf Hom}_{\hsolop}(sk_n\ycal'_+\wedge S_s^1, X_+)$. Then, using the
exactness
of sequences (\ref{due}), we conclude that ${\sf p}^*$ is bijective from
$$
{\sf Hom}_{\hsolop}(\ycal'_+,X_+)={\sf Hom}_{\hsolo}(\ycal',X_+)=
{\sf Hom}_{\hsolo}(\ycal, X_+),
$$
thus from ${\sf Hom}_{\hsolo}(\ycal, X)$. Theorem \ref{fibering} is completely
proved.
\enddemo

\medbreak

\begin{lemma}\label{solo}
Let $X\in\mathcal S_T$  and ${\sf p}:U\times\aff\to U$ be the
projection. Then the map
$$
{\sf p}^*:{\sf Hom}_{\fasc}(U,X)\to
{\sf Hom}_{\fasc}(U\times\aff,X)
$$
is bijective for every smooth scheme $U$ if and only if
\begin{equation}
{\sf p}^*_{k(u)}:{\sf Hom}_{\fasc}(Spec~k(u),X)\to
{\sf Hom}_{\fasc}(\mathbb{A}^1_{k(u)}, X)
\end{equation}
is, for every finite field extension $Spec~L\to Spec~k$,
${\sf p}_{L}:\mathbb{A}^1_{L}\to Spec~L$ being the projection.
\end{lemma}
\demo
We have just to prove that the bijectivity of $\sf p_L^*$ for
every $L$ finite extension of $k$ implies the bijectivity of $\sf p^*$
for every smooth scheme $U$. The morphism ${\sf p}:U\times\aff\to
U$ is a faithfully flat covering, thus, by faithfully flat descent we
have the following exact sequence of sets
$$
 0\to {\sf Hom}_{}(U,X)\stackrel{{\sf p}^*}{\to}
 {\sf Hom}_{}(U\times\aff,X)
\underset{p_2^*}{\stackrel{{\sf p}_1^*}{\rightrightarrows}}\\ {\sf
Hom}_{} ((U\times\aff)\times_U(U\times\aff),X).$$
In
order to prove the surjectivity of ${\sf p}^*$, we have to show that ${\sf
p}_1^*={\sf p}_2^*$. Notice that
$$
(U\times\aff)\times_U(U\times\aff)=U\times\mathbb{A}^2
$$
and ${\sf p}^*_1$ and ${\sf p}^*_2$ are induced by the projections on the
factors of $\mathbb{A}^2$ to $\aff$. Thus, given $\alpha\in {\sf
Hom}_{}(U\times\aff, X)$, we prove that
$$ \alpha\circ{\sf p}_1=\alpha\circ {\sf
p}_2: U\times\mathbb{A}^2\to X. $$
By hypothesis, any map $\aff_{L}\to X$ factors through
$Spec~L$ for any finite extension $L/k$. In particular,
$\alpha\circ{\sf p}_1$ and $\alpha\circ {\sf p}_2$ coincide on the
closed points of $U\times\mathbb{A}^2$. Since the union of all closed
points of $U\times\mathbb{A}^2$ is an everywhere dense subset for the
Zariski topology, we conclude that  $\alpha\circ{\sf p}_1=\alpha\circ
{\sf p}_2$.
\enddemo

\begin{corollary}\label{equivalenti}
Let  $X$ be a compact complex space. Then $X$ is Kobayashi hyperbolic if
and only if it is hyperbolic according to the definition \ref{spazio-iper}.
\end{corollary}

\demo Consequence of Theorem \ref{fibering}, Lemma
\ref{solo} and Brody's Theorem.
\enddemo
\begin{corollary}\label{brodyeq}
Let $X$ be a complex space, $C$ a closed complex subspace of $X$. Then $X$
is hyperbolic modulo $C$ in the sense of Brody if and only if $X/C$ is a hyperbolic  sheaf according to the
definition \ref{spazio-iper}.
\end{corollary}
\demo Let $\mathcal S_T$ be the site of complex spaces. By definition,
the sheaf of $\scal_T$ given by $Y\rightsquigarrow
\homs_{\fasc}(Y,X/C)$ is the associated sheaf for the strong
  topology to the presheaf which associates to a complex space $Y$ the
   colimit of

$$
\xymatrix{{\sf Hom}_{\mathcal S}(Y,
  C)\ar[r]\ar[d] & {\sf Hom}_{\mathcal S}(Y,X)\\{\sf Hom}_{\mathcal S}(Y,
  {\sf pt}).&}
$$
If  $X/C$ is a hyperbolic sheaf, then, by Theorem \ref{fibering}, we
obtain that the morphism
$$ {\sf Hom}_{\fasc}({\sf pt},
X/C)\to {\sf Hom}_{\fasc}(\compl, X/C)
$$
is a bijection. Assume, by
contradiction, that there exists a non constant holomorphic map $f:\compl
\to X$ such that $f(\compl)\not\subset C$. Then $f$ represents an element in
${\sf Hom}_{\fasc}(\compl, X/C)$ which is not in the image of ${\sf
Hom}_{\fasc}({\sf pt}, X/C)$ which is absurd. Conversely, if  $X$ is
Brody-hyperbolic modulo $C$, one has
$$ {\sf Hom}_{\mathcal S}(\compl,X)=(X-C)\amalg {\sf Hom}_{\mathcal S}
(\compl, C).
$$
On the other hand, we observe that ${\sf Hom}_{\fasc}(\compl, X/C)$ is
  precisely equal to the colimit of
$$
\xymatrix{{\sf Hom}_{\mathcal S}(\compl,
  C)\ar[r]\ar[d] & {\sf Hom}_{\mathcal S}(\compl,X)\\{\sf
    Hom}_{\mathcal S}(\compl, {\sf pt}).&}
$$
This follows from the fact that the new sections that we would get by
taking the associated sheaf are of the form $(f,g)$ where $f:U\to X$,
$g:V\to X$ are holomorphic maps, $\{U,\; V\}$ is an open covering of
$\compl$ (we may assume both $U$ and $V$ to be connected)
and $f(U\cap V)$, $g(U\cap V)$ are contained in $C$. In this
situation, we have that both $f(U)$ and $f(V)$ are contained in $C$, as
well. Therefore, $(f,g)=(U\to{\sf
  pt}, V\to{\sf pt})=\compl\to {\sf pt}$ and  we already
have this section in ${\sf Hom}_{\fasc}(\compl, X/C)$.
Consequently,
$$
{\sf Hom}_{\fasc}(\compl, X/C)=(X-C)\amalg {\sf
Hom}_{\mathcal S}(\compl,{\sf pt})
$$
but the latter set is ${\sf Hom}_{\fasc}({\sf pt},X/C)$, thus $X/C$ is
a hyperbolic sheaf by Theorem \ref{fibering}.
\enddemo
\medbreak \noindent Let us discuss some examples of hyperbolic resolutions of
complex spaces. Roughly speaking, $\mathfrak{Ip}(X)$ "enlarges" $X$ by
adding a simplicial structure which trivializes passing from $\hsolo$
to $\holo$. If $X$ is a Brody hyperbolic complex space, $\ip(X)$ is
isomorphic to $X$ in the category $\hsolo$. If  $X$ is not Brody
hyperbolic the simplicial structures added to $\ip(X)$ have the task
to "make constant" (up to simplicial homotopy, hence in $\hsolo$) all
morphisms $\compl\to X$. Passing from $\hsolo$ to $\holo$, $X$ and $\ip(X)$
become isomorphic objects .
\begin{esempi}\label{esem}
{\rm \begin{enumerate}
\item[1)] $\ip(\compl)$ is a simplicial sheaf isomorphic to a point in $\holo$.
Indeed, $\compl\cong {\sf pt}$ in $\holo$ and the hyperbolic resolutions
preserve affine equivalences. This fact is not surprising because if we
want to make all morphisms $\compl\to\compl$ homotopically constant, in
particular this must be true for the identity $\compl\to\compl$.
\item[2)]For the same reason, $\ip(\compl^n)\cong {\sf pt}$ in $\holo$
  for every $n\in{\mathbb N}$.
\item[3)] More generally, if $p:V\to X$ is a vector bundle,
%%%%%%%%%%%%parte nuova
  $\ip(p):\ip(V)\stackrel{\cong}{\to}\ip(X)$ in $\hsolo$ because $p$ is
  a $\compl$ weak equivalence.
Therefore, if $X$ is a hyperbolic complex space, then
$\ip(V)\cong X$ in $\hsolo$ and hence in $\holo$.
%%%%%%%%%%%%%%%%%
\item[4)] If $X$ is a complex space and $\ip(X)$ is represented by a
hyperbolic complex space $Y$, then $Y$ is unique up to isomorphisms
(cfr. the lemma below). In general, this is not the case; e.g. in the
next section we will show that $\ip(\mathbb{P}^n)$ cannot be
$\compl$-equivalent to a hyperbolic complex space.
\end{enumerate}}
\end{esempi}

In the case $\ip(\xcal)$ admits a hyperbolic complex space as
representative, then such a space is unique up to biholomorphism:
\begin{lemma}\label{unico}
Let $\xcal$ be a simplicial sheaf,  $Y$, $Y'$ hyperbolic complex
spaces such that
$$
\ip(\xcal)=[Y]_\holo=[Y']_\holo.
$$
Then $Y'$ and $Y$ are isomorphic complex spaces.
\end{lemma}
\demo Let $\mathcal S$ be the category of complex spaces. By hypothesis,
there exists an isomorphism $\phi: Y\cong Y'$ in $\holo$, namely a morphism
$\psi:Y'\to Y$ in $\holo$ such that $\psi\circ\phi={\sf id}_Y$ and
$\phi\circ\psi={\sf id}_{Y'}$ in $\holo$. Since  $Y$ e $Y'$ are complex
hyperbolic spaces, and in particular $\compl$-fibrant objects by Corollary
\ref{equivalenti} (see also the end of  Section \ref{base}), $\phi$ and
$\psi$ can be represented by morphisms $\phi':Y\to Y'$ and $\psi':Y'\to Y$ in
$\Delta^{\sf op}\fasc$. More precisely, we may suppose that $\phi'$ and
$\psi'$ are holomorphic maps, $Y$, $Y'$ being complex spaces and $\mathcal
S\hookrightarrow \Delta^{\sf op}\fasc$ being a full immersion.
Moreover, the fact
that $\phi$, $\psi$ are inverse to each other means that $\psi'\circ\phi'\sim
{\sf id}_Y$, $\phi'\circ \psi'\sim {\sf id}_{Y'}$ as holomorphic maps,
where $f\sim g$ if and only if there exists a holomorphic map
$H:W\times\compl\to V$ such that $H|_{W\times 0}=f$ e $H|_{W\times 1}=g$
(cfr. equation (\ref{quoz})). Since both $Y$, $Y'$ are hyperbolic, $H$ must
be constant along the fibres which are isomorphic to $\compl$, thus
$f\sim g$ if and only if $f=g$ as maps. In particular,
$\psi'\circ\phi'={\sf id}_Y$ and $\phi'\circ \psi'={\sf id}_{Y'}$.
\enddemo

%%%%%%%%%%%%%%%%%%%%%%%%%%%%%%%%%%%%%%%%%%%%%%%%%%%%%%%%%%%%%%%%%%%%%%%

In some cases, we can extend some results known for hyperbolic complex
spaces to hyperbolic sheaves:
\begin{lemma}
Let $\fasc$ be the category of sheaves of sets on the site of complex
spaces with the strong topology and $F$ be a hyperbolic sheaf. Then
$$
\homs_{\fasc}(\mathbb{P}^n,F)=F(\sf pt)
$$
for any $n\geq 1$. In other words, any sheaf map from $\mathbb{P}^n$ to a hyperbolic sheaf $F$
must be constant.
\end{lemma}
\demo
Consider the case $n=1$ first. Let $\pro=U_0\cup U_1$ be an open
covering with $U_0=\pro\setminus \{0\}$ and $U_1=\pro\setminus\{\infty\}$. Then the square
\begin{equation}\label{coca}
\begin{split}
\xymatrix{U_0\cap U_1\ar@{^(->}[r]^{\hskip 10pt i_0}\ar@{^(->}[d]^{i_1} &
  U_0\ar[d]\\ U_1\ar[r] & \pro}
\end{split}
\end{equation}
is cocartesian in the category of sheaves. Thus
\begin{equation}
\begin{split}
\homs_{\fasc}(\pro, F)=\lim \begin{pmatrix}
\xymatrix{\homs_{\fasc}(U_0\cap U_1, F)&
  \homs_{\fasc}(U_0,F)\ar[l]^{\hskip 20pt i_0^*}\\
  \homs_{\fasc}(U_1,F)\ar[u]^{i_1^*} &}
\end{pmatrix}.
\end{split}
\end{equation}
Since $U_0\cong U_1\cong\compl$, we have that
$$
\homs_{\fasc}(U_j,F)=\homs_{\fasc}({\sf pt}, F)=F({\sf pt})
$$
for $j=0,1$ because
of the theorem \ref{fibering}. Moreover,  $i_j^*$ are injective
because they have a retraction given by $f^*$ where $f:{\sf pt}\to U_0\cap
U_1$ is any point. We conclude the statement of the lemma in the case
of $\pro$ by noticing that the image of $i_0^*$ coincides with the
one of $i_1^*$. Consider now the open covering of $\mathbb{P}^n$
given by $U_0=\mathbb{P}^n\setminus\mathbb{P}^{n-1}\cong\compl^n$ and
$U_1=\mathbb{P}^n\setminus\{\infty\}$, where $\infty$ coincides with the point
$(0,0,\cdots, 0)\in U_0=\compl^n$. We get a cocartesian square like
(\ref{coca}) with $\mathbb{P}^n$ replacing $\pro$. The previous
argument carries through in the general case. The only thing to check
is that $\homs_{\fasc}(U_1,F)=F(\sf pt)$. Notice that the canonical
projection $p:U_1\to\mathbb{P}^{n-1}$ is a rank one vector
bundle. Locally on $\mathbb{P}^{n-1}$ (for the strong topology) it is
$V\times\compl$, where $V$ is an open affine of
$\mathbb{P}^{n-1}$. Hence
\begin{equation}
p^*_V:\homs_{\fasc}(V, F)\to \homs_{\fasc}(V\times\compl, F)
\end{equation}
are bijections for all $V$, since $F$ is hyperbolic. Glueing these
data for $V$ ranging on an open affine covering of $\mathbb{P}^{n-1}$,
we get that
\begin{equation}
p^*:\homs_{\fasc}(\mathbb{P}^{n-1},F)\to \homs_{\fasc}(U_1, F)
\end{equation}
is a bijection. By inductive assumption, we conclude that
$$
\homs_{\fasc}(U_1, F)=F(\sf pt).
$$
\enddemo

%%%%%%%%%%%%%%%%%%%%%%%%%%%%%%%%%%%%%%%%%%%%%%%%%%%%%%%%%%%%%%%%%%%%%%%%%

\section{Holotopy groups}\label{Holotopy/G}
Throughout this section, $S_T$ will denote the site of complex
spaces endowed with
the strong topology. A simplicial object of $S_T$ is, by definition, a {\it
simplicial complex space}. If we forget the complex structure, we could
study the objects of $S_T$  by means of the classical homotopy
groups. Isomorphism classes of homotopy
groups are invariant under homeomorphisms hence, a fortiori, under
biholomorphisms, however, they do not reflect the existence and the
properties of the complex
structure. A rather natural modification of the
definition of homotopy enables us to attach to every simplicial sheaf on
$S_T$ two families $\{\pi_{i,j}^{\sf par}(\xcal)\}_{i,j}$,
$\{\pi_{n,m}^{\sf iper}(z_1,z_2)(\xcal)\}_{m,n}$ of sets (cfr. Definition
\ref{gruppi})  which, for positive simplicial degrees, have a
canonical group structure and are invariant under biholomorphisms. We
will use these groups in Section \ref{appl} to show that there
exist complex spaces (e.g $\mathbb{P}^n$) whose
hyperbolic resolutions (cfr. Definition \ref{ris-iper}) are not isomorphic to
the class of hyperbolic complex spaces, not even in the category $\holo$.

Define the {\em parabolic circle} by $$ S^1_{\sf par}=\compl/(0\amalg 1),
$$ and we denote by $S^n_{\sf par}$ the sheaf $S^1_{\sf
par}\wedge\stackrel{n}{\cdots} \wedge S^1_{\sf par}$.

Let $D\subset\compl$ be the unit disc and $z_1\neq z_2$ two points of $D$.
We define the {\em hyperbolic circle} $S^1_{\sf iper}(z_1,z_2)$ by $$
S^1_{\sf iper}(z_1,z_2)=D/(z_1\amalg z_2) $$ and we denote by $S^n_{\sf
iper}(z_1,z_2)$ the sheaf $S^1_{\sf
iper}(z_1,z_2)\wedge\stackrel{n}{\cdots}\wedge S^1_{\sf iper}(z_1,z_2)$.

The quotients defining parabolic and hyperbolic circles are taken in the
category $\fasc$, even though, in view of a theorem of Cartan
(cfr. \cite{car}) the set theoretic quotients have a complex structure.
\begin{definition}\label{gruppi}
Let $\xcal$ be a simplicial sheaf on $S_T$. Define

\begin{equation}
\pi_{i,j}^{\sf par}(\xcal,x)={\sf Hom}_{\holop}((\mathbb{C}-0)^{\wedge
  j}\wedge S^{i-j}_{\sf par}, (\xcal,x))
\end{equation}
for $i\geq j\geq 0$,
\begin{equation}
\pi_{n,m}^{\sf iper}(z_1,z_2)(\xcal,x)={\sf Hom}_{\holop}(S_{\sf
  iper}^n(z_1,z_2)\wedge S_{\sf par}^m, (\xcal,x))
\end{equation}
for $n,m\geq 0$.

These sets are called respectively {\em parabolic holotopy pointed sets} of
$\xcal$ (or groups
in the case they are ) and {\em hyperbolic holotopy pointed sets} of
$\xcal$ (or groups in the case they are).
\end{definition}

\begin{oss}
{\em The definitions above are compatible with the classical ones of
algebraic topology. More precisely, let $\holo^{top}$ be the (unstable)
homotopy category of topological spaces (i.e. the localization of the
category of topological spaces with respect to the usual weak
equivalences); then we
have $\pi_n(X,x)={\sf Hom}_{\holo^{top}}((S^n,p),(X,x))$ for every
topological space $X$. Moreover, the topological realization functor
(cfr. Section \ref{trf}) provides functorial group homomorphisms
$\pi_{i,j}^{\sf par}(X,x)\to\pi_{i-j}(X,x)$ and $\pi_{n,m}^{\sf
  iper}(z_1,z_2)(\xcal,x)\to\pi_m(X,x)$ for any complex space $X$.
}
\end{oss}
\begin{lemma}\label{gruppo}
The sets $\pi_{i,j}^{\sf par}$, $\pi_{n,m}^{\sf iper}$ have a canonical
group structure for $i>j>0$ and $m>0$.
\end{lemma}
\demo The first step consists in proving that $S^1_{\sf par}\cong S^1_s$ in
$\holo$.

Consider the cofibration sequence
\begin{equation}
0\amalg 1\hookrightarrow \compl\to S^1_{\sf par}\to S^1_s\to\compl\wedge
S_s^1\to\cdots
\end{equation}
where $0\amalg 1$ and $\compl$ are pointed by $0$.  Since $\compl\cong
\pts$ in $\holo$, we have $\compl\wedge S_s^1\cong\pts$ in
$\holo$. Applying the functor ${\sf Hom}_{\holop}(~~, \zcal)$, in
view of Proposition \ref{mod-aff}, we obtain long exact
sequences of sets and, from these, the
isomorphism $$ {\sf Hom}_{\holop}(S^1_{\sf par},\zcal)\cong{\sf
Hom}_{\holop}(S_s^1,\zcal) $$ for every $\zcal\in\Delta^{\sf
op}_\bullet\fasc$. It follows that $S^1_{\sf par}\cong S^1_s$ in $\holo$.
The simplicial object $S^1_s$ is a cogroup (object) in $\hsolo$ (and
consequently in $\holo$). It is sufficient to observe that, if $a_{\sf
  str}$ is the associated sheaf for the strong topology, $S_s^1\cong
a_{\sf str}({\sf Sing}(S^1))$ in $\hsolo$ and that $S^1$ is a cogroup in
$\holo^{top}$ with projection

$$
{\sf p}:S^1\to S^1/(\{i\}\amalg\{-i\})\stackrel{\sf homeo}{\cong}S^1\vee S^1
$$
as structural map. Then, applying to ${\sf p}$ the functor
$a_{\sf str}({\sf Sing}(~~))$ we get a morphism

$$
[S_s^1]\to[S_s^1\vee S_s^1]=[S^1_s]\vee [S_s^1]
$$
in $\hsolo$ which satisfies the properties making it a
comultiplication. These properties are formulated in such a way to
induce on the sets ${\sf Hom}_{\hsolo}(S^1_s, \zcal)$ a natural group
structure. The same holds for ${\sf Hom}_{\holo}(S_s^1,\zcal)$.
\enddemo

%%%%%%%%%%%%%%%%%%%%%%%%%%%%%%%%%%%%%%%%%%%%%%%%%%%%%%%%%%%%%%%%%%%
%%%%%%%%%%%%%%%%%%%%%%%%%%%%%%%%%%%%%%%%%%%%%%%%%%%%%%%%%%%%%%%%%%

\begin{theorem}\label{iper-zero}
Let  $X$ be a hyperbolic sheaf. Then the groups
$\pi_{i,j}^{\sf par}(X,x)$, $\pi_{n,m}^{\sf iper}(X,x)$
vanish for $i-j>0$ and any $m>0$.
\end{theorem}

\demo We begin with proving that
 ${\sf Hom}_{\holop}(Y\wedge S^1_{\sf par}, X)=0$ for
every pointed complex space $(Y,\{y\})$. By definition (cfr. Section \ref{boh1}),
$$ Y\wedge S^1_{\sf par}=Y\times\compl/R $$
where $R$ is the complex space $Y\times (0\amalg 1)\cup y\times\compl$.
Since $Y\times\compl/R$ is a sheaf and $X$ is a fibrant space, by Theorem
\ref{fibering} we conclude that
\begin{multline}
{\sf Hom}_{\holo}(Y\times\compl/R,X)=
{\sf Hom}_{\fasc}(Y\times\compl/R,X)=\\
=\{f\in {\sf Hom}_{\fasc}
(Y\times\compl,X):f|_R=\text{\sf constant}\}.
\end{multline}

Moreover, since $Y\times\compl$ and $X$ are complex spaces, we have
$$ {\sf
Hom}_{\fasc}(Y\times\compl,X)={\sf Hom}_{\sf olom}(Y\times\compl,X).
$$ $X$
is Brody hyperbolic hence, for every $y\in Y$, the restriction
of a holomorphic map $f:Y\times\compl\to X$ to $y\times\compl$ is
constant. Furthermore, if $f\in {\sf Hom}_{\fasc}(Y\times\compl/R,X)$
then $f$ is constant on $Y\times 0\subset R$ and consequenty on the whole
$Y\times\compl$. It follows that, if $f$ is pointed, then $f$ must be constant
with image $x$, the base point of $X$. This shows that ${\sf
  Hom}_{\holop}(Y\wedge S^1_s, X)=x$ for any pointed complex space $Y$
and any hyperbolic pointed sheaf $X$.

We would like now to prove the same result with $Y$ being replaced by
a quotient sheaf $W=Y/Z$. Consider the following commutative diagram

\begin{equation}\label{grande}
\begin{split}
\xymatrix{ & R \ar[r]\ar[d] & {\sf pt}\ar[d] \\
Z\times\compl \ar[r]\ar[d] & Y\times\compl \ar[r] \ar[d]
&\dfrac{Y\times\compl}{R}\\
\compl\ar[r] & \dfrac{Y}{Z}\times\compl &}
\end{split}
\end{equation}
where the two squares are cocartesian. Consider now the two new
cocartesian squares

\begin{equation}
\begin{split}
\xymatrix{Z\times\compl\ar[r]\ar[d] & \dfrac{Y\times\compl}{R}\ar[d]
  &&
R\ar[r]\ar[d] & \sf{pt}\ar[d]\\
\compl\ar[r] & P && \dfrac{Y}{Z}\times\compl\ar[r] &
  \dfrac{W\times\compl}{R}=W\wedge S_{\sf par}^1.}
\end{split}
\end{equation}
By chasing the diagram (\ref{grande}) and using that
$\dfrac{Y}{Z}\times\compl$ and
$\dfrac{Y\times\compl}{R}$ are colimits of the relevant diagrams,
we find two sheaf maps $P\to W\wedge S_{\sf par}^1$ and $W\wedge S_{\sf
  par}^1\to P$ that are mutually inverses. By definition of $P$ we have
\begin{multline}
{\sf Hom}_{\holop}(W\wedge S_{\sf par}^1,X)=
{\sf  Hom}_{\holop}(P, X)={\sf Hom}_{\fasc_\bullet}(P,X)=\\
={\sf
  Hom}_{\fasc_\bullet}((Y\times\compl)/R,X)\times_{{\sf Hom}_
{\fasc_\bullet}(Z\times\compl,X)}{\sf Hom}_{\fasc_\bullet}(\compl, X).
\end{multline}
Recall that $(Y\times\compl)/R=Y\wedge S_{\sf par}^1$, thus, by the
first part of the proof of the proposition,  ${\sf
  Hom}_{\fasc_\bullet}((Y\times\compl)/R, X)=x$, the base point of
$X$. The same holds for ${\sf Hom}_{\fasc_\bullet}(\compl, X)$ because,
by assumption, $X$ is Brody hyperbolic. Therefore,
$$
{\sf Hom}_{\holop}(W\wedge S_{\sf par}^1,X)=x
$$
for any quotient sheaf
$W$, and in particular for
$$
W=(\compl\setminus 0)^{\wedge j}\wedge
S^{i-j-1}_{\sf par},\>\>W=S^n_{\sf iper}(z_1,z_2)\wedge S^{m-1}_{\sf par}
$$
(see Definition \ref{gruppi}).
\enddemo

\begin{corollary}\label{non}
Let $\xcal$ be a simplicial sheaf. Assume that $\pi_{i,j}^{\sf
par}(\xcal,x)\neq 0$ for $i-j>0$ or $\pi_{n,m}^{\sf
iper}(z_1,z_2)(\xcal,x)\neq 0$ for $m>0$. Then $\mathfrak{Ip}(\xcal)$ is
not $\compl$ weakly equivalent to a hyperbolic sheaf.
In particular, if $X$ is a complex space such that $\pi_{i,j}^{\sf
par}(\xcal,x)\neq 0$ for $i-j>0$ or $\pi_{n,m}^{\sf
iper}(z_1,z_2)(\xcal,x)\neq 0$ for $m>0$, then $X$  is
not a Brody hyperbolic complex space.
\end{corollary}
\demo The proofs for the two cases are similar so we consider only the case of
the parabolic holotopy groups. By definition, $$ \pi_{i,j}^{\sf
par}(\xcal,x)={\sf Hom}_{\holop}((\compl\setminus 0)^{\wedge j}\wedge
S_s^{i-j},(\xcal,x))$$ and this set is a quotient of $$ {\sf
Hom}_{\Delta^{\sf op}_\bullet\fasc}((\compl\setminus 0)^{\wedge j}\wedge
S^{i-j}_s,(\widetilde{\xcal},\widetilde{x})) $$ where
$(\widetilde{\xcal},\widetilde{x})$ is an $\compl$-fibrant pointed simplicial
sheaf $\compl$-weakly equivalent to $(\xcal,x)$. In
particular, we may assume that $\widetilde{\xcal}$ is the hyperbolic
resolution $\mathfrak{Ip}(\xcal)$ of $\xcal$. If $\mathfrak{Ip}(\xcal)$ were
$\compl$-weakly equivalent to a Brody hyperbolic complex space $X'$, then
$\pi_{i,j}^{\sf par}(\xcal,x)$ would be a quotient of
$$
{\sf Hom}_{\Delta^{\sf
op}_\bullet\fasc} ((\compl\setminus 0)^{\wedge j}\wedge S^{i-j}_s, (X',x'))
$$
which for $i-j>0$ consists only in the constant map with value $x$ (cfr.
 Theorem \ref{iper-zero}).
\enddemo

\begin{oss}
{\em As mentioned in section \ref{base}, to relate holotopy groups of a
complex space $X$ with
morphisms in $\Delta^{op}\fasc$ it is necessary to replace $X$ with its
hyperbolic model $\ip(X)$. Then we know that $\pi_{i,j}(X,x)$ will be a
quotient of the set $\homs_{\Delta^{op}_\bullet\fasc}(\mathbf{S}^{i,j},
\ip(X))$, where $\mathbf{S}^{i,j}$ is a pointed model of the relevant
sphere.}
\end{oss}

%%%%%%%%%%%%%%%%%%%%%%%%%%%%%%%%%%%%%%%%%%%%%%%%%%%%%%%%%%%%%%%%%%%%%%%
%%%%%%%%%%%%%%%%%%%%%%%%%%%%%%%%%%%%%%%%%%%%%%%%%%%%%%%%%%%%%%%%%%%%%%%
%%%%%%%%%%%%%%%%%%%%%%%%%%%%%%%%%%%%%%%%%%%%%%%%%%%%%%%%%%%%%%%%%%%%%%%
%%%%%%%%%%%%%%%%%%%%%%%%%%%%%%%%%%%%%%%%%%%%%%%%%%%%%%%%%%%%%%%%%%%%%%%

\section{The topological realization functor}\label{trf}

From now on, $\mathbb{CP}^n$ will denote the complex projective 
space seen as topological space.
We would like to compare objects in $\holo$ and $\holo(k)$ with the
topological spaces,
objects of the topological (unstable) homotopy category $\holo^{top}$.
We will show that there exists a functor $t^{olo}:\holo\to \holo^{top}$
which extends the functor which associates the underlying topological
space to a complex space. In the algebraic case extends the
corresponding functor which associates to an algebraic
variety over $\compl$, the topological space of its (Zariski) closed
points. The general case only applies to the site of smooth
varieties over a field $k$ which admits
an embedding $i$ in $\compl$. It involves passing from a simplicial sheaf
over $k$ to a simplicial sheaf over $\compl$ by means of $i^*$ (or,
more precisely, by means of its total left derived functor).
Recall that for a sheaf $F$ and a morphism of sites
$\phi:\mathfrak{S}_1\to \mathfrak{S}_2$, the sheaf
$\phi^*F$ on $\mathfrak{S}_1$ is defined as the associated sheaf to
the presheaf whose sections are
$(\phi^*F)(U)={\sf colim}_V F(V)$, where the colimit is taken over all
the morphisms $U\to \phi^{-1}V$ for $U\in\mathfrak{S}_1$ and any
$V\in\mathfrak{S}_2$.

\begin{definition}\label{real-topo}
Let $(\mathfrak{S},I)$ be a site with interval (cfr. Section 2.3
\cite{morvoe}) equipped with a realization functor
$r:\fs\to Top$ to the category of topological
spaces. Denote by $\holo(\fs)$ the $I$ homotopy category whose objects
are simplicial sheaves over $\fs$. Then a functor
$t_r:\holo(\fs)\to \holo^{top}$  with values in
the unstable homotopy category of topological spaces
is called a {\em topological realization
  functor} if the following properties are satisfied:
\begin{enumerate}
\item\label{one} if $X\in\Delta^{op}\mathcal{F}(\fs)$ is a
simplicial set, then
  the class  $t_r(X)$ can be represented by the geometric realization $|X|$;
\item\label{two} if $F$ is the sheaf $\homs_{\fs}(~~,X)$, where
$X\in\fs$, then $t_r(F)$ can be represented by $r(X)$;
\item\label{three} $t_r$ commutes with direct products and homotopy colimits.
\end{enumerate}
\end{definition}

\begin{theorem}\label{realiz}
The sites with interval $(\co, \compl)$ and $((Sm/k)_T,\aff_k)$
admit a topological realization functor, provided that $k$ can be
embedded in $\compl$ and $T$ is not finer then the flat topology.
\end{theorem}\vskip-.4cm
\demo
Let $\phi:(\fs_1,I_1)\to(\fs_2,I_2)$ be a reasonable countinuous map
of sites with interval (cfr. Definition 1.49 \cite{morvoe}).
Consider the functor $\phi^*:\Delta^{op}\mathcal{F}(\fs_2)\to
\Delta^{op}\mathcal{F}(\fs_1)$ obtained by applying the inverse image
functor on each component of the simplicial sheaf on $\fs_2$.
A classical result in model
categories assures the existence of the total left derivative between
homotopy categories of a functor, provided that such a functor
sends weak equivalences between cofibrant objects to weak equivalences.
In the case of $\phi^*$, we will not be able to prove this for every
simplicial sheaf on $\fs_2$ and the relevant $I$ model categories. However,
we can get the same result in the following way. We consider the
full category of $I_2$ local objects $\holo_{s,I_2}\subset\holo_s$
introduced in the
Theorem \ref{localizz}, which is equivalent to the $I$ homotopy
category  $\holo(\fs_2,I_2)$ by the same
theorem. Such a category has the property that a morphism is an $I_2$
weak equivalence if and only if it is a simplicial weak
equivalence. Thus, to show that $\phi^*$ admits a total left derived
functor between the $I$ homotopy categories, it is sufficient to
show that $\phi^*$ sends simplicial weak equivalences between
$I_2$ local objects (since every object is cofibrant) to simplicial
weak equivalences. Actually, since the property for a simplicial sheaf
$\xcal$ to be $I_2$ local is invariant under simplicial weak
equivalences on $\xcal$ (cfr. Definition \ref{a1loc})
to validate the same conclusion it suffices to show a weaker condition:
there exists  a (simplicial) resolution functor $\Phi$ and
a natural transformation $\Phi\to \sf id$ with the property that $\phi^*$
sends simplicial weak equivalences between simplicial shaves of the
kind $\Phi(\xcal)\to\Phi(\ycal)$  to simplicial weak equivalences
for all $I_2$ local simplicial sheaves $\xcal$ and $\ycal$.
But this is precisely the statement of Proposition 1.57.2. of
\cite{morvoe} where $\Phi$ is taken to be $\Phi_\Sigma$ introduced in
Lemma \ref{palle}. This shows the existence of the total left derived
functor $\mathbb{L}\phi^*:\holo(\fs_2,I_2)\to \holo(\fs_1,I_1)$ of
$\phi^*$. Explicitely, it is defined as follows: let $\xcal$ be a
simplicial sheaf over $\fs_2$, then $\mathbb{L}\phi^*(\xcal)$
is represented by the simplicial sheaf $\phi^*(\Phi_\Sigma({\mathfrak
 {Ip}}(\xcal)))$, where ${\mathfrak{Ip}}(\xcal)$ is the $I_2$ local
simplicial sheaf mentioned in the Theorem \ref{localizz}. This
definition is well
posed on $\holo(\fs_2,I_2)$ because of the above remarks and the fact that, if
$\xcal$ and $\xcal'$ represent the same class in $\holo(\fs_2,I_2)$, then
${\mathfrak{Ip}}(\xcal)$ and ${\mathfrak{Ip}}(\xcal')$ are simplicially
weak equivalent.

We will now consider the case of the site with interval
$(\co,\compl)$ since the algebraic case when $k=\compl$ is
entirely similar. We set the  realization functor $r:\fs\to Top$ to be
the one which associates the underlying topological space $X^{top}$
to a complex space $X$. Let $\pts$ be the site with interval
whose only nonempty object is the final object
$\sf pt$ and $\psi$ be the trivial morphism of sites with interval
$\pts\to \co$. Notice that a simplicial sheaf on $\pts$ is just
a simplicial set. We take the interval $I$ in $\pts$ to be the constant
simplicial set $\homs(\sf pt,\compl)$. The functor $\psi^*$ sends a
simplicial sheaf $\xcal$ on $\co$ to the simplicial set
$\xcal(\sf pt)$. Thus, $\psi^*(\compl)=I$ so that, in
particular, it is $I$ contractible. Because of this, the functor $\psi$
is said to be a {\em reasonable} continuous map of sites with
interval (cfr. Definition 3.16, \cite{morvoe}) and
$\mathbb{L}\psi^*$ has a particularly nice description:
$\mathbb{L}\psi^*(\xcal)$ is represented by the simplicial sheaf
$\psi^*\Phi_{\Sigma}(\xcal)$ where $\Sigma$ is the class of
representable sheaves on $\co$ (see Lemma 3.15. of \cite{morvoe}).
Let $\tlc$ be the category of locally contractible topological spaces.
We now endow the images of $\mathbb{L}\psi^*$ by a structure of
topological spaces in order to obtain a functor
$\holo(\Delta^{op}\fasci)\to\holo(\Delta^{op}\tlc)$.
If $\ycal$ is a simplicial sheaf that in each degree is a disjoint union of
representable sheaves $\amalg_{j\in J} Y_j$, then we set
$\theta\ycal:=\ycal^{top}$, where $\ycal^{top}$ is the simplicial
topological space having the topological space $\amalg_{j\in J}\ycal_j^{top}$ in
the corresponding degree. Since $\psi^*$ is reasonable,
$\mathbb{L}\psi^*(\xcal)=[\psi^*\ycal]$ where $\ycal$ is any
representable simplicial sheaf equipped with a simplicial weak
equivalence $\ycal\to \xcal$. Any two such models will give rise to
simplicially weak equivalent inverse images by Proposition
1.57.2. \cite{morvoe}, thus, in particular, $I$ weak equivalent. This
shows that the definition of $\theta$ induces a functor
$$
\holo=\holo(\Delta^{op}\fasci)\to\holo(\Delta^{op}\tlc)
$$
which we will call $\theta$, as well.

\begin{oss}
{\em $\holo(\Delta^{op}\tlc)$ is a full subcategory of
  $\holo(\Delta^{op}\boldsymbol{\mathcal{F}}_{open}(Tlc_{open}))$. The
latter category is the $I$ homotopy category taking as interval the sheaf
  $I=\homs_{cont}(~~,\compl)$. Such an interval is an object of
  $\Delta^{op}\tlc$, thus we can see $\holo(\Delta^{op}\tlc)$ as the
  localized category with respect to the $I=\compl$ weak equivalences,
  considering $\compl$ as constant simplicial topological space and no
  longer only as constant simplicial set.}
\end{oss}

\begin{proposition}\label{equ}
There is an equivalence of categories
$\gamma:\holo(\Delta^{op}\tlc)\cong\holo^{top}$.
\end{proposition}

\demo(sketch) It is a particular case of Proposition 3.3 of
\cite{morvoe}. Here we write the definition of the functor
$\gamma:\holo(\Delta^{op}\tlc)\to \holo^{top}$ which gives the equivalence of
categories. Let $\xcal$ be a
simplicial locally contractible topological space.
Since for any topological space $Z$ in
$\tlc$ there is an open covering $\amalg_i U_i\to Z$, with $U_i$
contractible for all $i$, by Lemma \ref{palle},
 $\xcal$ admits a (simplicial) weak equivalence
$\widetilde{\xcal}\to \xcal$ with
$\widetilde{\xcal}_j=\amalg_{i_j}U_{i_j}$. In turn,
$\tilde{\xcal}$ is $I$ weakly equivalent to $\xcal'$, where
$\xcal'$ is the {\em simplicial set} with $\xcal'_j=\amalg_{i_j}{\sf pt}$,
because $\widetilde{\xcal}$ and $\xcal'$
are termwise weakly equivalent and of Proposition 2.14, \cite{morvoe}.
The equivalence of categories is defined as
$[\xcal]\rightsquigarrow[|\xcal'|]$ where $|\xcal'|$ is the geometric
realization of the simplicial set $\xcal'$.
\enddemo

\begin{oss}\label{adopo}
{\em If $X$ is a topological space in $Tlc_{open}$, then $|\xcal'|$ is weakly
  equivalent to $|Sing_\bullet(X)|$. But this topological space is
  weakly equivalent to $X$ itself, thus the constant simplicial
  topological space $X$ is sent by $\gamma$ to a topological space
  weakly equivalent to $X$ in the classical sense of homotopy theory.}
\end{oss}

Let ${\mathbf D}$ be a small category and $\Delta^{op}\fasc^{\mathbf D}$ be the category of functors from ${\mathbf D}$ to $\Delta^{op}\fasc$. We will denote by 
$\rm{hocolim}({\mathbf D})$ a {\it homotopy colimit} of ${\mathbf D}$ on the category $\Delta^{op}\fasc$. That is a pair $(k,a)$ consisting in a functor $k:\Delta^{op}\fasc^{\mathbf D}\to\Delta^{op}\fasc$ which takes objectwise weak equivalences in $\Delta^{op}\fasc^{\mathbf D}$ to weak equivalences in $\Delta^{op}\fasc$ and a natural transformation $a:k\to \rm{colim}_{\mathbf D}$.
Such functor can be obtained by first taking a suitable cofibrant diagram replacement of an element in $\Delta^{op}\fasc^{\mathbf D}$ and composed with the ordinary colimit functor (cfr. \cite{bk}).

We set $t^{olo}:\holo\to\holo^{top}$ to be the functor $\gamma\circ\theta$.
Property (\ref{one}) of Definition \ref{real-topo} follows by
definition of $\gamma$. Property (\ref{two}) is a consequence of
Remark \ref{adopo}. As for the property (\ref{three}), we have that
$\psi^*$ commutes with limits by definition.
Since direct products in
the homotopy categories are represented by direct products of objects,
we have that $t^{olo}$ commutes with direct products. $\psi^*$ has a right
adjoint, namely $\psi_*$, thus it is right exact. Moreover, $\psi^*$
sends cofibrations (sectionwise injections) to cofibrations.
On the other hand, the same holds for the resolution functor of Lemma
\ref{palle} $\Phi$: if $i$ is a sectionwise injection, then $\Phi(i)$ is a
sectionwise injection by definition of $\Phi$; furthermore, $\Phi$
commutes with colimits, since its value on objects has been defined as a
colimit. In particular, if $\bf D$ is
a cofibrant diagram in $\Delta^{op}\fasci$, $\Phi(\bf D)$ is cofibrant and
we conclude that $\psi^*\Phi(\bf D)$ is cofibrant as well and also that
$\mathrm{colim}(\psi^*\Phi(\bf D))\cong \psi^*\Phi^*(\mathrm{colim}(\bf D))$.
This shows that, for any diagram $\bf D$,
$\mathbb{L}\psi^*(\mathrm{hocolim}(\bf D))=\mathrm{hocolim}(\mathbb{L}\psi^*(\bf D))$,
since the former class can be represented by
$\mathrm{colim}(\psi^*\Phi(\bf D'))$ for any cofibrant replacement
$\bf D'\stackrel{\sim}{\to} \bf D$ because $\psi^*\Phi(\bf D')$ is a cofibrant
diagram.

Therefore, $\theta$ commutes with homotopy colimits. Recall
that the equivalence $\gamma$ is defined to be the functor that, to a
class represented by a simplicial topological space $\xcal$, associates
the class in $\holo^{top}$ represented by
$|({\Phi_S(\xcal)})^{\sim}|$ where $S$ is the class of contractible
topological spaces and the operation $\sim$ replaces each contractible
topological space with a point. Because of the definition of $\Phi_S$
we see that $\sim$ sends injections to injections and commutes with
colimits. Before proceeding to
investigate the properties of the functor $|~~|$, we need to
recall the model structures involved in the categories. The functor
$|~~|$ is defined on the category of simplicial sets and takes values
in the category of topological spaces. The model structure for the
category of simplicial sets is:
 let $f:X\to Y$ be a map of simplicial sets, then $f$ is:\vskip-.7cm
{\em
\begin{enumerate}
\item a weak equivalence if $|f|$ is a weak homotopy equivalence (see below);
\item a cofibration if it is an injection;
\item a fibration if $f$ has the right lifting property with respect
  to acyclic cofibrations.
\end{enumerate}}
Let $X_0\hookrightarrow X_1\hookrightarrow X_2\to\cdots$ be a sequential
direct system of topological spaces such that for each $n$,
$(X_n, X_{n+1})$ is a relative CW complex. Then we will say that the
canonical function $X_0\hookrightarrow\mathrm{colim}X_i$ is a {\em
  generalized relative CW inclusion}. A continuous function between
topological spaces $f:X\to Y$ is
{\em
\begin{enumerate}
\item a weak equivalence if $f_*:\pi_*(X,x)\to \pi_*(Y,f(x))$ is a group
  isomorphism for $*\geq 1$ and a bijection of pointed sets if $*=0$;
\item a cofibration if it is a retract of a generalized relative
  CW inclusion;
\item a fibration if it is a Serre fibration.
\end{enumerate}}
The functor $|~~|$ preserves cofibrations and also it commutes
with colimits, because it has a right adjoint, namely the functor
${\sf Sing}(~~)$. In conclusion, the functor $\gamma$ commutes with homotopy
colimits, and so does the topological realization functor $t^{olo}$.

\enddemo

%%%%%%%%%%%%%%%%%%%%%%%%%%%%%%%%%%%%%%%%%%%%%%%%%%%%%%%%%%%%%%%%%%%
%%%%%%%%%%%%%%%%%%%%%%%%%%%%%%%%%%%%%%%%%%%%%%%%%%%%%%%%%%%%%%%%%%%

\subsection{Remarks on homotopy colimits}

The practical use of the topological realization functor requires few
remarks on the differences between (homotopy) colimits of diagrams in
the category $\holo^{top}$ and the category $\holo$. Let us consider
the colimit of the diagram
\begin{equation}\label{es1}
\begin{split}
\xymatrix{\compl\setminus 0\ar@{^(->}[r]\ar[d] &\compl\\ \sf pt.&}
\end{split}
\end{equation}
In the category of complex spaces, this is just a point. However, we
have previously  inferred  in this
manuscript that the colimit of such a diagram in the category of
sheaves on $\co$ is not (weakly equivalent to) the constant sheaf to a
point. Indeed, its class in the respective homotopy categories plays
the role of the two dimensional sphere $S^2=\mathbb C\pro$, or, more precisely,
of the sheaf represented by $\mathbb C\pro$, whose class is by no
means isomorphic to the one of the point. As a diagram of topological
spaces, its colimit is not a point, but it is not an appropriate
model for $S^2$. We should point out that the diagram (\ref{es1}) is a
cofibrant diagram for the affine model structure in the category
$\Delta^{op}\fasci$, but it is not cofibrant in the category of
topological spaces for the model structure defined above.
This apparent oddness disappears if we consider
{\em homotopy} colimits instead. For instance,
$$
t^{olo}(S^n_{\mathsf{par}})\cong t^{olo}(S^n_{\mathsf{iper}}(z_1,z_2))\cong
S^n
$$
for any $n\geq 0$ and $z_i\in D$, $t^{olo}(\compl\setminus 0)\cong S^1$ and we have natural
maps of pointed sets (respectively of groups)
$$
\pi_{i,j}^{\mathsf{par}}(\xcal,x)\to
\pi_i(t^{olo}(\xcal),t^{olo}(x))
$$
and
$$
\pi_{n,m}^{\mathsf{iper}}
(\xcal,x)\to\pi_{n+m}(t(\xcal),t^{olo}(x)).
$$

%%%%%%%%%%%%%%%%%%%%%%%%%%%%%%%%%%%%%%%%%%%%%%%%%%%%%%%%%%%%%%%%%%%%%%
%%%%%%%%%%%%%%%%%%%%%%%%%%%%%%%%%%%%%%%%%%%%%%%%%%%%%%%%%%%%%%%%%%%%%%
%%%%%%%%%%%%%%%%%%%%%%%%%%%%%%%%%%%%%%%%%%%%%%%%%%%%%%%%%%%%%%%%%%%%%%

\section{Some applications}\label{appl}

%%%%%%%%%%%%%%%%%%%%%%%%%%%%%%%%%%%%%%%%%%%%%%%%%
%%%%%%%%%%%%%%%%%%%%%%%%%%%%%parte nuova%%%%%%%%%
%%questa sezione ha subito varie aggiunte%%%%%%%%
%%%%%%%%%%%%%%%%%%%%%%%%%%%%%%%%%%%%%%%%%%%%%%%%%
In this last section we are going to consider few applications of the
theory developed so far. We will begin with examples of complex spaces that
are not $\compl$ weakly equivalent to any complex hyperbolic space.

\begin{definition}
We will say that a complex space is {\em weakly hyperbolic} if is
$\compl$ weakly equivalent to a Brody hyperbolic complex space.
\end{definition}

We recall a preliminary result (cfr. Lemma 2.15 \cite{morvoe}):
\begin{lemma}\label{sfera1}
The pointed simplicial sheaf $(\compl\setminus 0)\wedge S^1_{\sf par}$ is
canonically weakly equivalent to $\pro$.
\end{lemma}
\demo
Consider the diagram $\bf D$
\begin{equation}
\begin{split}
\xymatrix{(\compl\setminus\{0\},\{1\})\ar[r]\ar[d] & (\compl,\{1\})\\
(\compl\setminus \{0\},\{1\})\wedge\Delta[1].& }
\end{split}
\end{equation}
If $\bf D'$ is another diagram
\begin{equation}
\begin{split}
\xymatrix{\mathcal X\ar[r]^f\ar[d]^i& Y\\
\mathcal Z&}
\end{split}
\end{equation}
in $\holo$ then $\rm{colim}_{\bf D}\cong\text\rm{colim}_{{\bf D}'}$ in
$\holo$ if there exists a morphism of diagrams ${\bf D}\to {\bf D}'$ such
that the morphisms are weak affine equivalences. Consider the diagrams ${\bf
D}'$ e ${\bf D}''$
\begin{equation}
\begin{split}
\xymatrix{(\compl\setminus\{0\},\{1\})\ar[r]\ar[d] &{\sf pt }&&
  (\compl\setminus\{0\},\{1\})\ar[r]\ar[d] & (\compl,\{1\})\\
(\compl\setminus\{0\},\{1\})\wedge \Delta[1] & && {\sf pt} &}
\end{split}
\end{equation}
and the morphisms $$ f=((\compl,\{1\})\to {\sf pt}, {\sf id}):{\bf D}\to {\bf
D}'\>\>\>\>\> g=({\sf id}, (\compl\setminus\{0\},\{1\})\wedge\Delta[1]\to {\sf pt}):{\bf
D}\to {\bf D}''. $$
 The morphisms $f$ and $g$ induce affine weak equivalences
 colim$_{\bf D}\to$\rm{colim}$_{{\bf D}'}$ and
colim$_{\bf D}\to$\rm{colim}$_{{\bf D}''}$. Identifying colim$_{{\bf
D}'}$ with $(\compl\setminus\{0\},\{1\})\wedge S_s^1$ and colim$_{{\bf D}''}$ with
$\compl/(\compl\setminus\{0\})$, we conclude that $$ (\compl\setminus\{0\},\{1\})\wedge S_s^1\cong\compl/(\compl\setminus\{0\}). $$ The square
\begin{equation}
\begin{split}
\xymatrix{\compl\setminus\{0\}\ar[r]\ar[d] & \compl\ar[d]\\
\pro\setminus\{\infty\}\ar[r] &\pro}
\end{split}
\end{equation}
is cocartesian in $\fasc$, hence the cofibres of horizontal morphisms are
isomorphic. We derive
  $$
  \compl/(\compl\setminus\{0\})\cong \pro/(\pro\setminus\{\infty\})
  $$
   in $\fasc$. But
  $$
  \pro/(\pro\setminus\{\infty\})\cong \pro
  $$
  in $\holo$, since $\pro\setminus\{\infty\}\cong {\sf pt}$ in $\holo$.
\begin{oss}\label{sfera2}
{\em In the proof of Lemma \ref{gruppo} we have already seen that
$S^1_{\mathsf{par}}$ is weakly equivalent to $S^1_s$.}
\end{oss}

%%%%%%%%%%%%%%%%%%%%%%%%%%%%%%%%%%%%%%%%%%%%%%%%%%%%%%%%%%%%%%%%%%%%%%%%

We are now going to apply the theory developed so far to prove that
\begin{theorem}\label{no}
For any $n>0$, $\mathbb{P}^n$ is not weakly hyperbolic. In other words,
$\ip(\mathbb{P}^n)$ cannot be represented in $\holo$ by
a Brody hyperbolic complex space.
\end{theorem}
\demo
In view of Corollary \ref{non}, it is sufficient to show that

$$ \pi_{2,1}^{\mathsf{par}}(\mathbb{P}^n,\infty)={\sf
Hom}_{\holo_\bullet}((\compl\setminus\{0\})\wedge S^{1}_{\sf par},
(\mathbb{P}^n,\{\infty\}))\neq 0
$$
or equivalently, by Lemma \ref{sfera1} and Remark \ref{sfera2}, that
$$
\homs_{\holo_\bullet}(\pro,(\mathbb{P}^n,\{\infty\}))\neq 0.
$$
Our candidate to represent a nonzero class is the canonical embedding
$i:\mathbb{P}^1\hookrightarrow \mathbb{P}^n$. The topological realization
yields a group homomorphism
$$
t:\pi_{2,1}^{\mathsf{par}}(\mathbb{P}^n,\infty)\to\pi_{2}(\mathbb{CP}^n,
\infty).
$$
$t^{olo}(i):\mathbb{CP}^1\hookrightarrow \mathbb{CP}^n$ is
the canonical inclusion and not null homotopic, since $\mathbb{CP}^n$
is obtained by $\mathbb{CP}^1$ by attaching cells of dimension $4$ and
above, hence it is an equivalence up to dimension $2$ and in
particular
$$
t^{olo}(i)_*:
\mathbb{Z}=\pi_2(\mathbb{CP}^1,\infty)\to\pi_2(\mathbb{CP}^n,\infty)
$$
is an isomorphism. In conclusion $t[i]\neq 0$, thus
$[i]\neq 0\in \pi_{2,1}^{\mathsf{par}}(\mathbb{P}^n,\infty)$.
\enddemo

%%%%%%%%%%%%%%%%%%%%%%%%%%%%%%%%%%%%%%%%%%%%%%%%%%%%%%%%%%%%%%%%%%%%%%
%%%%%parte nuova

\begin{proposition}\label{riv}
Let $X$ be a complex space and $p:\widetilde{X}\to X$ a connected covering
complex space. Assume that $X$ is weakly hyperbolic and let $f:\compl\to X$
be a nonconstant holomorphic function. Then for any lifting $\tilde{f}$ of
$f$ to $\widetilde{X}$, $\tilde{f}(\compl)$ contains just one point in each
fiber of $p$  or equivalently $p|_{\tilde{f}(\compl)}$ is a biholomorphism
for any such $f$ and $\tilde{f}$.
\end{proposition}
\demo
Let $X$ be weakly hyperbolic. Assume, by a contradiction, that there exist a nonconstant holomorphic function
$f:\compl\to X$ and a lifting $\tilde{f}:\compl\to \tilde{X}$ such that $a\neq b\in p^{-1}(x)$, $x\in X$,  $a,b\in \tilde{f}(\compl)$. For
the purposes of this proof, we
can assume that $\tilde{f}(0)=a$ and $\tilde{f}(1)=b$.
Then we have the following commutative diagram:
\begin{equation}
\begin{split}
\xymatrix{\compl\ar[r]^{\tilde{f}}\ar[d]^q & \tilde{X}\ar[d]^p\\
\compl/\{0\}\amalg \{1\}\ar@{.>}[r]^{\>\>\>\>\alpha} & X}
\end{split}
\end{equation}
where $\alpha$ sends the class of $\{0\}\amalg \{1\}$ to $x\in X$. We
have that $[\alpha]\neq 0\in \pi_{1,0}^{\mathsf{par}}(X,x)$. Indeed,
$[\alpha^{top}]\neq 0\in \pi_1(X^{top},x)$. Consider the composition
$$
[0,1]\stackrel{g}{\to}\compl/\{0\}\amalg \{1\}\stackrel{\alpha^{top}}{\to}
X^{top},
$$
where $g$ is a path from $0$ to $1$ in $\compl$. If
$\alpha^{top}\circ g$ is not homotopic to a constant relatively to
$\{0,1\}$, then $\alpha^{top}$ is not homotopic to a constant. But, by
construction, $\alpha^{top}\circ g$ lifts uniquely to a path in
$\widetilde{X}^{top}$ starting from $a$ and ending in $b$, hence
$\alpha^{top}\circ g$ cannot be homotopic to a
constant relatively to $\{0,1\}$. This shows that $\pi_1(X^{top},x)\neq 0$ which is absurd since $X$ is weak hyperbolic.
\enddemo

The Proposition \ref{riv} in particular implies the following
\begin{corollary}
 Any complex space $X$ whose
universal covering space is $\compl^n$ for some $n\geq 1$, is not
weakly hyperbolic.
\end{corollary}
\demo Let $p:\mathbb C^n\to X$ be the universal covering of $X$. Let $a\neq
b\in p^{-1}(x)$, $x\in X$. A complex line $l\subset\mathbb C^n$ passing
through $a,b$ provides a homorphic map $f:\mathbb C\to X$ which does not
satisfy the conclusion of Proposition \ref{riv}.
\enddemo
%%%%%%%%%%%%%%%%%%%%%%%%%%%%%%%%%%%%%%%%%%%%%%%%%%%%%%%%%%%%%%%%%%%%%%%
%%%%%%%%%%%%%%%%%%%%%%%%%%%%%%%%%%%%%%%%%%%%%%%%%%%%%%%%%%%%%%%%%%%%%%%

{\sc Universit\`a degli Studi di Milano - Bicocca - Piazza dell'Ateneo

Nuovo, 1 - 20126, Milano, Italy}

 e-mail: { \tt simone.borghesi@unimib.it}

{\sc Scuola Normale Superiore di Pisa, Piazza dei Cavalieri

7, 56126 Pisa, Italy}

e-mail: {\tt g.tomassini@sns.it}

\end{document}